\documentclass[11pt]{amsart}

\date{\today}

\usepackage[margin=1in]{geometry}
\usepackage{setspace}
\usepackage[hang,flushmargin,symbol*]{footmisc}
\usepackage{amsmath}
\usepackage{amsthm}
\usepackage{amssymb}
\usepackage{mathtools}
\usepackage{enumitem}
\usepackage{graphicx}
\usepackage{caption}
\usepackage[labelformat=simple,labelfont={}]{subcaption}
\usepackage{color}
\definecolor{darkblue}{rgb}{0, 0, .6}
\usepackage{hyperref}
\hypersetup{
	colorlinks=true,
	linkcolor=darkblue,
	anchorcolor=darkblue,
	citecolor=darkblue,
	pagecolor=darkblue,
	urlcolor=darkblue,
	pdftitle=darkblue,
	pdfauthor=darkblue
}

\newtheorem{theorem}{Theorem}[subsection]
\newtheorem{lemma}[theorem]{Lemma}

\newtheorem{corollary}[theorem]{Corollary}
\newtheorem{proposition}[theorem]{Proposition}

\theoremstyle{definition} 
\newtheorem{definition}[theorem]{Definition}
\newtheorem{example}[theorem]{Example}
\newtheorem{remark}[theorem]{Remark}

\newcommand{\Z}{\mathbb{Z}}

\newcommand{\C}{\widetilde{C}}
\renewcommand{\O}{\mathcal{O}}
\newcommand{\E}{\mathcal{E}}

\newcommand{\wtri}{\vartriangle}
\newcommand{\btri}{\blacktriangle}
\renewcommand{\a}{\mathbf{a}}
\DeclareMathOperator{\TL}{TL}
\DeclareMathOperator{\DTL}{\mathbb{D}TL}
\renewcommand{\P}{\mathcal{P}}
\newcommand{\V}{\mathcal{V}}
\newcommand{\D}{\mathbb{D}}

\newcommand{\Diag}{\mathcal{D}}
\newcommand{\wcirc}{\circ}

\newcommand{\bcirc}{\bullet}

\renewcommand{\H}{\mathcal{H}}
\DeclareMathOperator{\FC}{FC}

\begin{document}

\title[Diagram calculus for a type affine $C$ Temperley--Lieb algebra, I]{Diagram calculus for a type affine $C$\\ Temperley--Lieb algebra, I}

\author{Dana C.~Ernst}
\address{
Department of Mathematics and Statistics,
Northern Arizona University PO Box 5717,
Flagstaff, AZ 86011-5717, USA
}
\email{Dana.Ernst@nau.edu}



\begin{abstract}

In this paper, we present an infinite dimensional associative diagram algebra that satisfies the relations of the generalized Temperley--Lieb algebra having a basis indexed by the fully commutative elements (in the sense of Stembridge) of the Coxeter group of type affine $C$.  Moreover, we provide an explicit description of a basis for the diagram algebra.  In the sequel to this paper, we show that this diagrammatic representation is faithful.  The results of this paper and its sequel will be used to construct a Jones-type trace on the Hecke algebra of type affine $C$, allowing us to non-recursively compute leading coefficients of certain Kazhdan--Lusztig polynomials.

\end{abstract}

\maketitle


\begin{section}{Introduction}\label{sec:intro}

The (type $A$) Temperley--Lieb algebra $\TL(A)$, invented by H.N.V.~Temperley and E.H.~Lieb in 1971~\cite{Temperley.H;Lieb.E:A}, is a finite dimensional associative algebra which first arose in the context of statistical mechanics.  R.~Penrose and L.H.~Kauffman showed that this algebra can be realized as a diagram algebra~\cite{Kauffman.L:B, Penrose.R:A}, that is, an associative algebra with a basis given by certain diagrams, in which the multiplication rule in the algebra is given by applying local combinatorial rules to the diagrams.  

In 1987, V.F.R.~Jones showed that $\TL(A)$ occurs naturally as a quotient of the type $A$ Hecke algebra~\cite{Jones.V:B}.  Given a Coxeter group $W$, the associated Hecke algebra has a basis indexed by the elements of $W$ and relations that deform the relations of $W$ by a parameter $q$. The realization of the Temperley--Lieb algebra as a Hecke algebra quotient was generalized by J.J.~Graham in~\cite{Graham.J:A} to the case of an arbitrary Coxeter system.  In Section~\ref{subsec:TL-algebras}, we define the generalized Temperley--Lieb algebra of type $\C_{n}$, denoted $\TL(\C_{n})$, in terms of generators and relations and describe a special basis, called the monomial basis, which is indexed by the fully commutative elements (defined in Section~\ref{subsec:FC}) of the underlying Coxeter group.

The goal of this paper is to introduce a diagrammatic representation of the Temperley--Lieb algebra (in the sense of Graham) of type $\C$.  The motivation behind this is that a realization of $\TL(\C_n)$ can be of great value when it comes to understanding the otherwise purely abstract algebraic structure of the algebra.  Moreover, studying these generalized Temperley--Lieb algebras often provides a gateway to understanding the Kazhdan--Lusztig theory of the associated Hecke algebra.  Loosely speaking, the generalized Temperley--Lieb algebra retains some of the relevant structure of the Hecke algebra, yet is small enough that computation of the leading coefficients of the notoriously difficult to compute Kazhdan--Lusztig polynomials is often much simpler.

In this paper, we construct an infinite dimensional associative diagram algebra $\D_n$ that satisfies the relations of $\TL(\C_n)$.  In Sections~\ref{sec:diagram algebras} and~\ref{sec:simple and admissible}, we establish our notation and introduce all of the necessary terminology required to define $\D_n$, and once this has been done it is trivial to verify that the relations of $\TL(\C_n)$ are satisfied and that there is a surjective algebra homomorphism from $\TL(\C_n)$ to $\D_n$ (Proposition~\ref{prop:surjective homomorphism}).  However, due to length considerations, the injectivity of the homomorphism is resolved in the sequel to this paper~\cite{Ernst.D:D}.

One of the major obstacles to proving that our diagrammatic representation is faithful is having a description of a basis for $\D_n$.  In Section~\ref{subsec:admissible}, we define the $\C$-admissible diagrams by providing a combinatorial description of the allowable edge configurations involving diagram decorations.  Our main result (Theorem~\ref{thm:module admissibles is D_n}) comes at the end of a sequence of technical lemmas and states that the $\C$-admissible diagrams form a basis for $\D_n$.  Finally, in Section~\ref{sec:closing}, we discuss the implications of our results and future research.

With the exception of type $\widetilde{A}$, all other generalized Temperley--Lieb algebras with known diagrammatic representations are finite dimensional.  In the finite dimensional case, counting arguments are employed to prove faithfulness, but these techniques are not available in the type $\C$ case since $\TL(\C_n)$ is infinite dimensional.  Instead, we will make use of the author's classification in~\cite{Ernst.D:B} of the non-cancellable elements in Coxeter groups of types $A$, $B$, and $\C$ (also see~\cite[Chapters 3--5]{Ernst.D:A}).  The classification of the non-cancellable elements in a Coxeter group of type $\C$ provides the foundation for inductive arguments used to prove the faithfulness of $\D_n$.  Once injectivity has been established, the diagram algebra introduced in this paper will be the first faithful representation of an infinite dimensional non-simply-laced generalized Temperley--Lieb algebra (in the sense of Graham).

This paper is an adaptation of the author's Ph.D. thesis, titled \textit{A diagrammatic representation of an affine $C$ Temperley--Lieb algebra}~\cite{Ernst.D:A}, which was directed by Richard M.~Green at the University of Colorado at Boulder.  However, the notation has been improved and some of the arguments have been streamlined.  In particular, the author's thesis describes a  framework for constructing a large class of diagram algebras and is more general than what often appears in the literature.  For the sake of length, we omit here the general construction and focus on our diagram algebra of interest.
\end{section}


\begin{section}{Preliminaries}


\begin{subsection}{Coxeter groups}\label{subsec:coxeter groups}

A \emph{Coxeter group} is a group $W$ with a distinguished set of
generating involutions $S$ having presentation
\[
\langle s_1,\dots,s_n\mid (s_is_j)^{m(s_{i},s_{j})}=1\rangle,
\]
where $m:S\times S\to \mathbb{N}$ is a function and $m(s_{i},s_{j})=1$ if and only if $i=j$.  It turns out that the elements of $S$ are distinct as group elements, and that $m(s, t)$ is the order of $st$.  Any minimum length expression for $w \in W$ in terms of the generators is called a \emph{reduced expression} (all reduced expressions for $w$ have the same length).  The pair $(W,S)$ is called a \emph{Coxeter system}.

Given a Coxeter system $(W,S)$, the associated \emph{Coxeter graph} $\Gamma$ is the graph with vertex set $S$ and edges $\{s,t\}$ for each $m(s,t)\geq 3$. Moreover, each edge is labeled with its corresponding $m$-value, although it is customary to omit the label if $m(s,t)=3$. Given a Coxeter graph $\Gamma$, we can uniquely reconstruct the corresponding Coxeter system $(W,S)$.  In this case, we say that the corresponding Coxeter system is of type $\Gamma$, and denote the Coxeter group and distinguished generating set by $W(\Gamma)$ and $S(\Gamma)$, respectively.

The main focus of this paper will be the Coxeter systems of types $B_n$ and $\C_n$, which are defined by the Coxeter graphs in Figures~\ref{Fig002} and~\ref{Fig003}, respectively, where $n\geq 2$.

\begin{figure}[!ht]
\centering
\subcaptionbox{Coxeter graph of type $B_n$.\label{Fig002}}{\includegraphics[scale=.95]{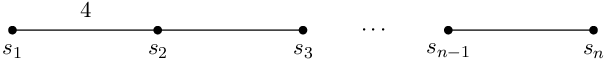}}
\vspace{1em}
\subcaptionbox{Coxeter graph of type $\C_n$.\label{Fig003}}{\includegraphics[scale=.95]{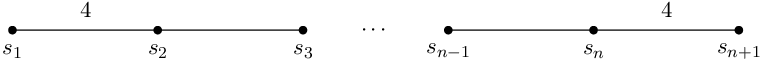}}
\caption{Coxeter graphs.}
\end{figure}

We can obtain $W(B_{n})$ from $W(\C_{n})$ by removing the generator $s_{n+1}$ and the corresponding relations~\cite[Chapter 5]{Humphreys.J:A}.  We also obtain a Coxeter group of type $B$ if we remove the generator $s_{1}$ and the corresponding relations.  To distinguish these two cases, we let $W(B_{n})$ denote the subgroup of $W(\C_{n})$ generated by $\{s_{1}, s_{2}, \dots, s_{n}\}$ and we let $W(B'_{n})$ denote the subgroup of $W(\C_{n})$ generated by $\{s_{2}, s_{3}, \dots, s_{n+1}\}$.  It is well-known that $W(\C_{n})$ is an infinite Coxeter group while $W(B_{n})$ and $W(B'_{n})$ are both finite~\cite[Chapters 2 and 6]{Humphreys.J:A}.

\end{subsection}


\begin{subsection}{Fully commutative elements}\label{subsec:FC}

Let $(W,S)$ be a Coxeter system of type $\Gamma$ and let $w \in W$. According to Stembridge~\cite{Stembridge.J:B}, $w$ is \emph{fully commutative} (FC) if and only if no reduced expression for $w$ contains a subword of the form $ststs \cdots$ of length $m(s,t) \geq 3$.  We will denote the set of all FC elements of $W$ by $\FC(W)$ or $\FC(\Gamma)$.

The elements of $\FC(\C_{n})$ are precisely those whose reduced expressions avoid subwords of the following types:
\begin{enumerate}
\item $s_{i}s_{j}s_{i}$ for $|i-j|=1$ and $1< i,j < n+1$;
\item $s_{i}s_{j}s_{i}s_{j}$ for $\{i,j\}=\{1,2\}$ or $\{n,n+1\}$.
\end{enumerate}
The FC elements of $W(B_{n})$ and $W(B'_{n})$ avoid the respective subwords above.

By~\cite[Theorem~5.1]{Stembridge.J:B}, $W(\C_{n})$ contains an infinite number of FC elements, while $W(B_{n})$ (and hence $W(B'_n)$) contains finitely many.  There are examples of infinite Coxeter groups that contain a finite number of FC elements (e.g., $W(E_n$) is infinite for $n\geq 9$, but contains only finitely many FC elements~\cite[Theorem~5.1]{Stembridge.J:B}).

\end{subsection}


\begin{subsection}{Generalized Temperley--Lieb algebras}\label{subsec:TL-algebras}

Given a Coxeter graph $\Gamma$, we can form an associative algebra, $\TL(\Gamma)$ (in the sense of Graham~\cite{Graham.J:A}), which we call the Temperley--Lieb algebra of type $\Gamma$.  For a complete description of the construction of $\TL(\Gamma)$, see~\cite{Ernst.D:A,Graham.J:A,Green.R:P}.  For our purposes it suffices to define $\TL(\C_{n})$ in terms of generators and relations.  We are using~\cite[Proposition~2.6]{Green.R:P} (also see~\cite[Proposition~9.5]{Graham.J:A}) as our definition. 

\begin{definition}\label{def:TL(C)}
The \emph{Temperley--Lieb algebra of type $\C_{n}$}, denoted $\TL(\C_{n})$, is the unital algebra generated by $\{b_{1}, b_{2}, \dots, b_{n+1}\}$ with defining relations
\begin{enumerate}
\item $b_{i}^{2}=\delta b_{i}$ for all $i$, where $\delta$ is an indeterminate;
\item $b_{i}b_{j}=b_{j}b_{i}$ if $|i-j|>1$;
\item $b_{i}b_{j}b_{i}=b_{i}$ if $|i-j|=1$ and $1< i,j < n+1$;
\item $b_{i}b_{j}b_{i}b_{j}=2b_{i}b_{j}$ if $\{i,j\}=\{1,2\}$ or $\{n,n+1\}$.
\end{enumerate}
In addition, $\TL(B_{n})$ (respectively, $\TL(B'_{n})$) is generated as a unital algebra by $\{b_{1}, b_{2}, \dots, b_{n}\}$ (respectively, $\{b_{2}, b_{3}, \dots, b_{n+1}\}$) with the relations above.
\end{definition}

It is known that we can consider $\TL(B_{n})$ and $\TL(B'_{n})$ as subalgebras of $\TL(\C_{n})$ in the obvious way.

Note that when $\TL(\C_{n})$ is considered as a quotient of the Hecke algebra of type $\C_{n}$ with indeterminate $v$, the indeterminate $\delta$ is defined to be the Laurent polynomial $v+v^{-1}$.

Let $s_{x_1}s_{x_2}\cdots s_{x_r}$ be a reduced expression for $w \in \FC(\C_{n})$, where each $x_i \in \{1,\ldots, n+1\}$.  Define the element $b_w \in\TL(\C_{n})$ via
\[
b_w=b_{s_{x_1}}b_{s_{x_2}}\cdots b_{s_{x_r}}.
\]
It is well-known (and follows from~\cite[Proposition~2.4]{Green.R:P}) that the set $\{b_{w}: w \in \FC(\C_{n})\}$ forms a $\Z[\delta]$-basis for $\TL(\C_{n})$.  This basis is referred to as the \emph{monomial basis} or ``$b$-basis.''  

If $(W,S)$ is a Coxeter system of type $\Gamma$, the associated Hecke algebra $\H(\Gamma)$ is an algebra with a basis indexed by the elements of $W$ and relations that deform the relations of $W$ by a parameter $q$. In general, $\TL(\Gamma)$ is a quotient of $\H(\Gamma)$, having several bases indexed by the FC elements of $W$~\cite[Theorem~6.2]{Graham.J:A}.  Except for in the case of type $A$, there are many Temperley--Lieb type quotients that appear in the literature.  That is,  some authors define a Temperley--Lieb algebra to be a different quotient of $\H(\Gamma)$ than the one we are interested in.  In particular, the blob algebra of~\cite{Martin.P;Saleur.H:A} is a smaller Temperley--Lieb type quotient of $\H(B_{n})$ than $\TL(B_{n})$.  Also, the symplectic blob algebra of~\cite{Green.R;Martin.P;Parker.A:A} and ~\cite{Martin.P;Green.R;Parker.A:A} is a finite rank quotient of $\H(\C_{n})$, whereas, $\TL(\C_{n})$ is of infinite rank.  Furthermore, despite being infinite dimensional, the two-boundary Temperley--Lieb algebra of~\cite{Gier.J;Nichols.A:A} is a quotient of $\H(\C_n)$ different from $\TL(C_{n})$.  Typically, authors that study these usually smaller Temperley--Lieb type quotients are interested in representation theory, whereas our motivation is Kazhdan--Lusztig theory.  

\end{subsection}

\end{section}


\begin{section}{Diagram algebras}\label{sec:diagram algebras}

The goal of this section is to familiarize the reader with the necessary background on diagram algebras.  It is important to note that there is currently no rigorous definition of the term ``diagram algebra.''  Our diagram algebras possess many of the same features as those already appearing in the literature, however the typical developments are too restrictive to accomplish the task of finding a faithful diagrammatic representation of the infinite dimensional Temperley--Lieb algebra (in the sense of Graham) of type $\C$.  Yet, our approach is modeled after~\cite{Green.R:M},~\cite{Green.R;Martin.P;Parker.A:A},~\cite{Jones.V:A}, and~\cite{Martin.P;Green.R;Parker.A:A}.


\begin{subsection}{Undecorated diagrams}

First, we discuss undecorated diagrams and their corresponding diagram algebras.

\begin{definition}\label{def:k-box}
Let $k$ be a nonnegative integer.  The \emph{standard $k$-box} is a rectangle with $2k$ marks points, called \emph{nodes} (or \emph{vertices}) labeled as in Figure~\ref{Fig073}.  We will refer to the top of the rectangle as the \emph{north face} and the bottom as the \emph{south face}.
\end{definition}

Sometimes, it will be useful for us to think of the standard $k$-box as being embedded in the plane.  In this case, we put the lower left corner of the rectangle at the origin such that each node $i$ (respectively, $i'$) is located at the point $(i,1)$ (respectively, $(i,0)$).

\begin{figure}[!ht]
\centering
\includegraphics[scale=.95]{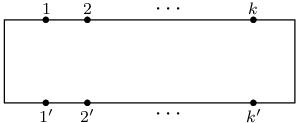}
\caption{The standard $k$-box.}\label{Fig073}
\end{figure}

The next definition summarizes the construction of the ordinary Temperley--Lieb pseudo diagrams.  

\begin{definition}\label{def:T_k(emptyset)}
A \emph{concrete pseudo $k$-diagram} consists of a finite number of disjoint curves (planar), called \emph{edges}, embedded in the standard $k$-box with the following restrictions.  The nodes of the box are the endpoints of edges, which meet the box transversely.  All other edges must be closed (isotopic to circles) and disjoint from the box.  We define an equivalence relation on the set of concrete pseudo $k$-diagrams.  Two concrete pseudo $k$-diagrams are \emph{(isotopically) equivalent} if one concrete diagram can be obtained from the other by isotopically deforming the edges such that any intermediate diagram is also a concrete pseudo $k$-diagram.  A \emph{pseudo $k$-diagram} (or an \emph{ordinary Temperley-Lieb pseudo-diagram}) is defined to be an equivalence class of equivalent concrete pseudo $k$-diagrams.  We denote the set of pseudo $k$-diagrams by $T_{k}(\emptyset)$.
\end{definition}

\begin{example}\label{ex:pseudo diagram}
The diagram in Figure~\ref{Fig074} is an example of a concrete pseudo 5-diagram.
\end{example}

\begin{figure}[!ht]
\centering
\includegraphics[scale=.95]{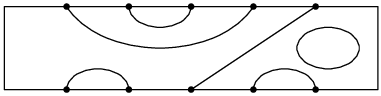}
\caption{A concrete pseudo 5-diagram.}\label{Fig074}
\end{figure}

\begin{remark}\label{rem:diagram repn}
When representing a pseudo $k$-diagram with a drawing, we pick an arbitrary concrete representative among a continuum of equivalent choices.  When no confusion can arise, we will not make a distinction between a concrete pseudo $k$-diagram and the equivalence class that it represents.  
\end{remark}

We will refer to a closed curve occurring in the pseudo $k$-diagram as a \emph{loop edge}, or simply a \emph{loop}.  The diagram in Figure~\ref{Fig074} has a single loop.  Note that we used the word ``pseudo'' in our definition to emphasize that we allow loops to appear in our diagrams.  Most examples of diagram algebras occurring in the literature ``scale away'' loops that appear.  There are loops in the diagram algebra that we are interested in preserving, so as to obtain infinitely many diagrams.  The presence of $\emptyset$ in the definition above is to emphasize that the edges of the diagrams are undecorated.  In the next section, we allow for the presence of decorations.

Let $d$ be a diagram.  If $d$ has an edge $e$ that joins node $i$ in the north face to node $j'$ in the south face, then $e$ is called a \emph{propagating edge from $i$ to $j'$}.  (Propagating edges are often referred to as ``through strings'' in the literature.)  If a propagating edge joins $i$ to $i'$, then we will call it a \emph{vertical propagating edge}.  If an edge is not propagating, loop edge or otherwise, it will be called \emph{non-propagating}.  

If a diagram $d$ has at least one propagating edge, then we say that $d$ is \emph{dammed}.  If, on the other hand, $d$ has no propagating edges (which can only happen if $k$ is even), then we say that $d$ is \emph{undammed}.  Note that the number of non-propagating edges in the north face of a diagram must be equal to the number of non-propagating edges in the south face.  We define the function $\a: T_{k}(\emptyset) \to \Z^{+}\cup \{0\}$ via
\[
\a(d)=\text{ number of non-propagating edges in the north face of } d.
\]
There is only one diagram with $\a$-value $0$ having no loops; namely the diagram $d_{e}$ that appears in Figure~\ref{Fig076}.  The maximum value that $\a(d)$ can take is $\lfloor k/2 \rfloor$.  In particular, if $k$ is even, then the maximum value that $\a(d)$ can take is $k/2$, i.e., $d$ is undammed.  On the other hand, if $\a(d)=\lfloor k/2 \rfloor$ while $k$ is odd, then $d$ has a unique propagating edge.

\begin{figure}[!ht]
\centering
\includegraphics[scale=.95]{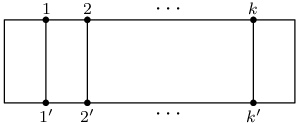}
\caption{The only diagram having $\a$-value 0 and no loops.}\label{Fig076}
\end{figure}

We wish to define an associative algebra that has the pseudo $k$-diagrams as a basis.

\begin{definition}\label{def:P_k(emptyset)}
Let $R$ be a commutative ring with $1$.  The associative algebra $\P_{k}(\emptyset)$ over $R$ is the free $R$-module having $T_{k}(\emptyset)$ as a basis, with multiplication defined as follows.  If $d, d' \in T_{k}(\emptyset)$, the product $d'd$ is the element of $T_{k}(\emptyset)$ obtained by placing $d'$ on top of $d$, so that node $i'$ of $d'$ coincides with node $i$ of $d$, rescaling vertically by a factor of $1/2$ and then applying the appropriate translation to recover a standard $k$-box.  (For a proof that this procedure does in fact define an associative algebra see~\cite[\textsection 2]{Green.R:M} and~\cite{Jones.V:A}.)
\end{definition}

We will refer to the multiplication of diagrams as \emph{diagram concatenation}.  The (ordinary) Temperley--Lieb diagram algebra (see~\cite{Green.R:N, Green.R:M, Jones.V:A, Penrose.R:A}) can be easily defined in terms of this formalism.

\begin{definition}\label{def:DTL(A_n)}
Let $\DTL(A_{n})$ be the associative $\Z[\delta]$-algebra equal to the quotient of $\P_{n+1}(\emptyset)$ by the relation depicted in Figure~\ref{Fig077}.
\end{definition}

\begin{figure}[!ht]
\centering
$\begin{tabular}[c]{@{} c@{}}
\includegraphics[scale=.95]{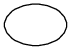} 
\end{tabular}
= \delta$
\caption{The defining relation of $\DTL(A_{n})$.}\label{Fig077}
\end{figure}

It is well-known that $\DTL(A_{n})$ is the free $\Z[\delta]$-module with basis given by the elements of $T_{n+1}(\emptyset)$ having no loops. The multiplication is inherited from the multiplication on $\P_{n+1}(\emptyset)$ except we multiply by a factor of $\delta$ for each resulting loop and then discard the loop.  We will refer to $\DTL(A_{n})$ as the \emph{(ordinary) Temperley--Lieb diagram algebra}.

\begin{example}
Figure~\ref{Fig078--Fig079} depicts the product of three basis diagrams of $\DTL(A_{4})$.

\begin{figure}[!ht]
$\begin{tabular}[c]{l}
\includegraphics[scale=.95]{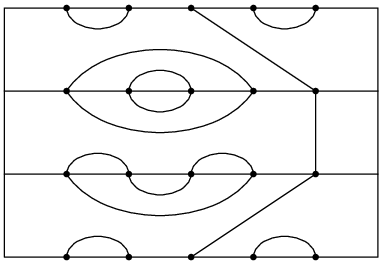}
\end{tabular}
= \  \delta^{3} \ \begin{tabular}[c]{l}
\includegraphics[scale=.95]{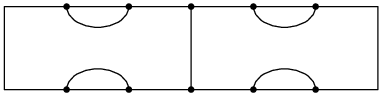}
\end{tabular}$
\caption{An example of multiplication in $\DTL(A_{4})$.}\label{Fig078--Fig079}
\end{figure}
\end{example}

As $\Z[\delta]$-algebras, the Temperley--Lieb algebra $\TL(A_{n})$ that was briefly discussed in Section~\ref{sec:intro} is isomorphic to $\DTL(A_{n})$.  Moreover, each loop-free diagram from $T_{n+1}(\emptyset)$ corresponds to a unique monomial basis element of $\TL(A_{n})$.  For more details, see~\cite{Kauffman.L:B} and~\cite{Penrose.R:A}.

\end{subsection}


\begin{subsection}{Decorated diagrams}

We wish to adorn the edges of a diagram with elements from an associative algebra having a basis containing $1$.  First, we need to develop some terminology and lay out a few restrictions on how we decorate our diagrams.

Let $\Omega=\{\bcirc, \btri, \wcirc, \wtri\}$ and consider the free monoid $\Omega^{*}$.  We will use the elements of $\Omega$ to adorn the edges of a diagram and we will refer to each element of $\Omega$ as a \emph{decoration}.  In particular, $\bcirc$ and $\btri$ are called \emph{closed decorations}, while $\wcirc$ and $\wtri$ are called \emph{open decorations}.   Let $\mathbf{b}=x_{1}x_{2}\cdots x_{r}$ be a finite sequence of decorations in $\Omega^{*}$.  We say that $x_{i}$ and $x_{j}$ are \emph{adjacent} in $\mathbf{b}$ if $|i-j|=1$ and we will refer to $\mathbf{b}$ as a \emph{block} of decorations of \emph{width} $r$.  Note that a block of width $1$ is just a single decoration.  The string $\bcirc~\bcirc~\btri~\wcirc~\bcirc~\wtri~\bcirc$ is an example of a block of width 7 from $\Omega^*$.

We have several restrictions for how we allow the edges of a diagram to be decorated, which we will now outline.  Let $d$ be a fixed concrete pseudo $k$-diagram and let $e$ be an edge of $d$.

\begin{enumerate}[label=\rm{(D0)}]

\item \label{D0} If $\a(d)=0$, then $e$ is undecorated.

\end{enumerate}

\noindent In particular, the unique diagram $d_{e}$ with $\a$-value 0 and no loops is undecorated.

Subject to some restrictions, if $\a(d)>0$, we may adorn $e$ with a finite sequence of blocks of decorations $\mathbf{b}_{1}, \dots, \mathbf{b}_{m}$ such that adjacency of blocks and decorations of each block is preserved as we travel along $e$.  

If $e$ is a non-loop edge, the convention we adopt is that the decorations of the block are placed so that we can read off the sequence of decorations as we traverse $e$ from $i$ to $j'$ if $e$ is propagating, or from $i$ to $j$ (respectively, $i'$ to $j'$) with $i < j$ (respectively, $i' < j'$) if $e$ is non-propagating.

If $e$ is a loop edge, reading the corresponding sequence of decorations depends on an arbitrary choice of starting point and direction round the loop. We say two sequences of blocks are \emph{loop equivalent} if one can be changed to the other or its opposite by any cyclic permutation. Note that loop equivalence is an equivalence relation on the set of sequences of blocks.  So, the sequence of blocks on a loop is only defined up to loop equivalence.  That is, if we adorn a loop edge with a sequence of blocks of decorations, we only require that adjacency be preserved.

Each decoration $x_{i}$ on $e$  has coordinates in the $xy$-plane.  In particular, each decoration has an associated $y$-value, which we will call its \emph{vertical position}.  

If $\a(d)\neq 0$, then we also require the following.

\begin{enumerate}[label=\rm{(D\arabic*)}]

\item \label{D1}  All decorated edges can be deformed so as to take closed decorations to the left wall of the diagram and open decorations to the right wall simultaneously without crossing any other edges.

\item \label{D2} If $e$ is non-propagating (loop edge or otherwise), then we allow adjacent blocks on $e$ to be conjoined to form larger blocks.

\item \label{D3}If $\a(d)>1$ and $e$ is propagating, then as in~\ref{D2}, we allow adjacent blocks on $e$ to be conjoined to form larger blocks.

\item \label{unusual} If $\a(d)=1$, then we have the following.

\begin{enumerate}
\item All decorations occurring on propagating edges must have vertical position lower (respectively, higher) than the vertical positions of decorations occurring on the (unique) non-propagating edge in the north face (respectively, south face) of $d$.

\item If a block on a propagating edge contains decorations occurring at vertical positions $y_{1}$ and $y_{2}$ with $y_{1} < y_{2}$, then no other propagating edge may contain decorations at vertical positions in the interval $(y_{1},y_{2})$.

\item Two adjacent blocks occurring on a propagating edge may be conjoined to form a larger block as long as (b) is not violated.
\end{enumerate}

\end{enumerate}

We call a block \emph{maximal} if its width cannot be increased by conjoining it with another block without violating~\ref{unusual}.

Requirement~\ref{D1} is related to the concept of ``exposed'' that appears in the context of the Temperley--Lieb algebra of type $B$~\cite{Green.R:N,Green.R:H,Green.R:M}.  The general idea is to mimic what happens in the type $B$ case on both the east and west ends of the diagrams.  Note that~\ref{unusual} is an unusual requirement for decorated diagrams.  We require this feature to ensure faithfulness of our diagrammatic representation on the monomial basis elements of $\TL(\C_{n})$ indexed by the type I elements of the Coxeter group of type $\C_{n}$ (see~\cite{Ernst.D:B}).

\begin{definition}
A \emph{concrete LR-decorated pseudo $k$-diagram} is any concrete $k$-diagram decorated by elements of $\Omega$ that satisfies conditions~\ref{D0}--\ref{unusual}.
\end{definition}

\begin{example}\label{ex:decorated diagrams}
Here are a few examples.
\begin{enumerate}[label=\rm{(\alph*)}]
\item \label{ex:decorated diagram1} The diagram in Figure~\ref{Fig080} is an example of a concrete LR-decorated pseudo $5$-diagram.  In this diagram, there are no restrictions on the relative vertical position of decorations since the $\a$-value is greater than 1.  The decorations on the unique propagating edge can be conjoined to form a maximal block of width 4.

\item \label{ex:decorated diagram2} The diagram in Figure~\ref{Fig081} is another example of a concrete LR-decorated pseudo $5$-diagram, but with $\a$-value 1.  We use the horizontal dotted lines to indicate that the three closed decorations on the leftmost propagating edge are in three distinct blocks.  We cannot conjoin these three decorations to form a single block because there are decorations on the rightmost propagating edge occupying vertical positions between them.  Similarly, the open decorations on the rightmost propagating edge form two distinct blocks that may not be conjoined.

\item \label{ex:decorated diagram3} Finally, the diagram in Figure~\ref{Fig082} is an example of a concrete LR-decorated pseudo $6$-diagram with maximal $\a$-value and no propagating edges.
\end{enumerate}
\end{example}

\begin{figure}[!ht]
\centering
\subcaptionbox{\label{Fig080}}{\includegraphics[scale=.95]{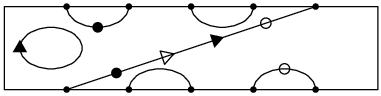}}
\qquad
\subcaptionbox{\label{Fig081}}{\includegraphics[scale=.95]{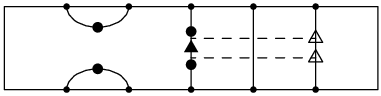}}\\
\subcaptionbox{\label{Fig082}}{\includegraphics[scale=.95]{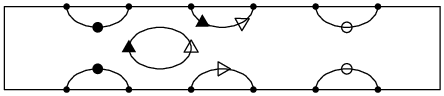}}
\caption{Examples of concrete LR-decorated pseudo diagrams.}
\end{figure}

Note that an isotopy of a concrete LR-decorated pseudo $k$-diagram $d$ that preserves the faces of the standard $k$-box may not preserve the relative vertical position of the decorations even if it is mapping $d$ to an equivalent diagram.  The only time equivalence is an issue is when $\a(d)=1$.  In this case, we wish to preserve the relative vertical position of the blocks.  We define two concrete pseudo LR-decorated $k$-diagrams to be \emph{$\Omega$-equivalent} if we can isotopically deform one diagram into the other such that any intermediate diagram is also a concrete LR-decorated pseudo $k$-diagram.  Note that we do allow decorations from the same maximal block to pass each other's vertical position (while maintaining adjacency).  

\begin{definition}\label{def:T_{k}^{LR}(Omega)}
An \emph{LR-decorated pseudo $k$-diagram} is defined to be an equivalence class of $\Omega$-equivalent concrete LR-decorated pseudo $k$-diagrams.  We denote the set of LR-decorated diagrams by $T_{k}^{LR}(\Omega)$.
\end{definition}

As in Remark~\ref{rem:diagram repn}, when representing an LR-decorated pseudo $k$-diagram with a drawing, we pick an arbitrary concrete representative among a continuum of equivalent choices.  When no confusion will arise, we will not make a distinction between a concrete LR-decorated pseudo $k$-diagram and the equivalence class that it represents. 

\begin{remark}\label{rem:LR-decorated}
We make several observations.
\begin{enumerate}
\item The set of LR-decorated diagrams $T_{k}^{LR}(\Omega)$ is infinite since there is no limit to the number of loops that may appear.

\item \label{rem:closed left open right} If $d$ is an undammed LR-decorated diagram, then all closed decorations occurring on an edge connecting nodes in the north face (respectively, south face) of $d$ must occur before all of the open decorations occurring on the same edge as we travel the edge from the left node to the right node. Otherwise, we would not be able to simultaneously deform decorated edges to the left and right.  Furthermore, if an edge joining nodes in the north face of $d$ is adorned with an open (respectively, closed) decoration, then no non-propagating edge occurring to the right (respectively, left) in the north face may be adorned with closed (respectively, open) decorations.  We have an analogous statement for non-propagating edges in the south face.

\item \label{rem:undammed dec loops} Loops can only be decorated by both types of decorations if $d$ is undammed.  Again, we would not be able to simultaneously deform decorated edges to the left and right, otherwise.

\item \label{rem:both decs one prop} If $d$ is a dammed LR-decorated diagram, then closed decorations (respectively, open decorations) only occur to the left (respectively, right) of and possibly on the leftmost (respectively, rightmost) propagating edge.  The only way a propagating edge can have decorations of both types is if there is a single propagating edge, which can only happen if $k$ is odd.
\end{enumerate}
\end{remark}

\begin{example}
The diagram of Figure~\ref{Fig082} is an example that illustrates conditions (\ref{rem:closed left open right}) and (\ref{rem:undammed dec loops}) of Remark~\ref{rem:LR-decorated}, while the diagram of Figure~\ref{Fig080} illustrates condition (\ref{rem:both decs one prop}).
\end{example}

\begin{definition}\label{def:P_{k}^{LR}(Omega)}
We define $\P_{k}^{LR}(\Omega)$ to be the free $\Z[\delta]$-module having the LR-decorated pseudo $k$-diagrams $T_{k}^{LR}(\Omega)$ as a basis.  
\end{definition}

We define multiplication in $\P_{k}^{LR}(\Omega)$ by defining multiplication in the case where $d$ and $d'$ are basis elements, and then extend bilinearly.  To calculate the product $d'd$, concatenate $d'$ and $d$ (as in Definition~\ref{def:P_k(emptyset)}).  While maintaining $\Omega$-equivalence, conjoin adjacent blocks.  We claim that the multiplication just defined turns $\P_{k}^{LR}(\Omega)$ into a well-defined associative $\Z[\delta]$-algebra.  To justify this claim we require the following lemma.

\begin{lemma}\label{lem:a-value=1}
Let $d$ be diagram with $\a(d)=1$.  Suppose that the unique non-propagating edge in the north face of $d$ joins $i$ to $i+1$.  Let $d'$ be any other diagram with $\a(d')>0$.  Then $\a(d'd)=1$ if and only if $\a(d')=1$ and the unique non-propagating edge in the south face of $d'$ joins either (a) $(i-1)'$ to $i'$, (b) $i'$ to $(i+1)'$, or (c) $(i+1)'$ to $(i+2)'$.
\end{lemma}

\begin{proof}
First, assume that $\a(d'd)=1$.  It is a general fact that $\a(d'd)\geq \a(d')$, which implies that $\a(d')=1$.  

Conversely, assume that $\a(d)=1$ and that the unique non-propagating edge in the south face of $d'$ joins either (a) $(i-1)'$ to $i'$, (b) $i'$ to $(i+1)'$, or (c) $(i+1)'$ to $(i+2)'$.  

Assume that we are in situation (a).  Suppose that the propagating edge leaving node $(i+1)'$ in the south face of $d'$ is connected to node $j$ in the north face.  Also, suppose that the propagating edge leaving node $i-1$ in the north face of $d$ is connected to node $l'$ in the south face.  Then $d'd$ has a propagating edge joining node $j$ to node $l'$.  Furthermore, the only non-propagating edge in the north (respectively, south) face of $d'd$ is the same as the unique non-propagating edge in the north (respectively, south) face of $d'$ (respectively, $d$).  It follows that $\a(d'd)=1$.  

Next, assume we are in case (b).  Then $d'd$ has one more loop than the sum total of loops from $d'$ and $d$.  Furthermore, the only non-propagating edge in the north (respectively, south) face of $d'd$ is the same as the unique non-propagating edge in the north (respectively, south) face of $d'$ (respectively, $d$), and so $\a(d'd)=1$.  

Finally, if we are in situation (c), then the proof that $\a(d'd)=1$ is symmetric to case (a).
\end{proof}

It is quickly seen that concatenating two diagrams that satisfy~\ref{D1} will result in a diagram that satisfies the same conditions.  The claim that $\P_{k}^{LR}(\Omega)$ is a well-defined associative $\Z[\delta]$-algebra now follows from arguments in~\cite[\textsection 3]{Martin.P;Green.R;Parker.A:A} and Lemma~\ref{lem:a-value=1} above.  The only case that requires serious consideration is when multiplying two diagrams that both have $\a$-value $1$.  If $\a(d)=\a(d')=1$ while $\a(d'd)>1$, then there are no concerns.  However, if $\a(d'd)=1$, then according to Lemma~\ref{lem:a-value=1}, if the unique non-propagating edge $e'$ in the south face of $d'$ joins $i'$ to $(i+1)'$, it must be the case that the unique non-propagating edge $e$ in the north face of $d$ joins either (a) $i-1$ to $i$, (b) $i$ to $i+1$, or (c) $i+1$ to $i+2$.  If (a) or (c) occurs, then the only blocks that get conjoined are the blocks on $e$ and $e'$, which presents no problems.  If (b) occurs, then we get a loop edge and we conjoin the blocks from $e$ and $e'$.  As a consequence, it is possible that the block occurring on a propagating edge of $d'$ having the lowest vertical position may be conjoined with the block occurring on a propagating edge of $d$ having the highest vertical position.  This can only happen if these two edges are joined in $d'd$, and regardless, presents no problems.

We remark that since the set of LR-decorated diagrams is infinite, $\P_{k}^{LR}(\Omega)$ is an infinite dimensional algebra.  

\end{subsection}


\begin{subsection}{Diagrammatic relations}

Our immediate goal is to define a quotient of $\P_{k}^{LR}(\Omega)$ having relations that are determined by applying local combinatorial rules to the diagrams. 

Let $R=\Z[\delta]$ and define the algebra $\V$ to be the quotient of $R\Omega^{*}$ by the following relations:
\begin{enumerate}
\item $\bcirc~\bcirc~=\btri$;
\item $\bcirc~\btri=\btri~\bcirc~=~2~\bcirc$;
\item $\wcirc~\wcirc~=~\wtri$;
\item $\wcirc\wtri~=~\wtri\wcirc~=~2~\wcirc$.
\end{enumerate}
The algebra $\V$ is associative and has a basis consisting of the identity and all finite alternating products of open and closed decorations.  

For example, in $\V$ we have
\[
\bcirc~\bcirc~\wcirc~\bcirc~\wcirc~\wcirc~\bcirc~=\btri~\wcirc~\bcirc~\wtri~\bcirc,
\]
where the expression on the right is a basis element, while the expression on the left is a block of width 7, but not a basis element.  We will refer to $\V$ as our \emph{decoration algebra}.

The point is that there is no interaction between open and closed symbols.  It turns out that if $\delta=1$, the algebra $\V$ is equal to the free product of two rank 3 Verlinde algebras.  For more details, see Chapter 7 of the author's Ph.D. thesis~\cite{Ernst.D:A}.

\begin{definition}\label{def:big diagram alg}
Let $\widehat{\P}_{k}^{LR}(\Omega)$ be the associative $\Z[\delta]$-algebra equal to the quotient of $\P_{k}^{LR}(\Omega)$ by the relations depicted in Figure~\ref{Fig083--Fig095}, where the decorations on the edges represent adjacent decorations of the same block.
\end{definition}

\begin{figure}[!ht]
\centering
\subcaptionbox{\label{Fig083--Fig084}}{$\begin{tabular}[c]{@{}c@{}}
\includegraphics[scale=.95]{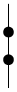}
\end{tabular}
= 
\begin{tabular}[c]{@{} c@{}}
\includegraphics[scale=.95]{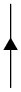}
\end{tabular}$}
\qquad
\subcaptionbox{\label{Fig085--Fig086}}{$\begin{tabular}[c]{@{}c@{}}
\includegraphics[scale=.95]{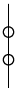}
\end{tabular}
= 
\begin{tabular}[c]{@{} c@{}}
\includegraphics[scale=.95]{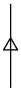}
\end{tabular}$}
\qquad
\subcaptionbox{\label{Fig087--Fig089}}{$\begin{tabular}[c]{@{} c@{}}
\includegraphics[scale=.95]{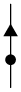}
\end{tabular}
=  
\begin{tabular}[c]{@{} c@{}}
\includegraphics[scale=.95]{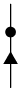}
\end{tabular}
=   2
\begin{tabular}[c]{@{} c@{}}
\includegraphics[scale=.95]{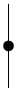}
\end{tabular}$}
\qquad
\subcaptionbox{\label{Fig090--Fig092}}{$\begin{tabular}[c]{@{} c@{}}
\includegraphics[scale=.95]{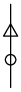}
\end{tabular}
=  
\begin{tabular}[c]{@{} c@{}}
\includegraphics[scale=.95]{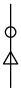}
\end{tabular}
=   2
\begin{tabular}[c]{@{} c@{}}
\includegraphics[scale=.95]{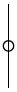}
\end{tabular}$}
\qquad
\subcaptionbox{\label{Fig093--Fig095}}{$\begin{tabular}[c]{@{} c@{}}
\includegraphics[scale=.95]{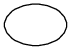}
\end{tabular}
=  
\begin{tabular}[c]{@{} c@{}}
\includegraphics[scale=.95]{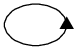}
\end{tabular}
=  
\begin{tabular}[c]{@{} c@{}}
\includegraphics[scale=.95]{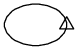}
\end{tabular}
=\   \delta$}
\caption{The defining relations of $\widehat{\P}_{k}^{LR}(\Omega)$.}\label{Fig083--Fig095}
\end{figure}

Note that with the exception of the relations involving loops, multiplication in $\widehat{\P}_{k}^{LR}(\Omega)$ is inherited from the relations of the decoration algebra $\V$.   Also, observe that all of the relations are local in the sense that a single reduction only involves a single edge.  As a consequence of the relations in Figure~\ref{Fig083--Fig095}, we also have the relations of Figure~\ref{Fig096--Fig099}.

\begin{figure}[!ht]
\centering
\subcaptionbox{\label{Fig096--Fig097}}{$\begin{tabular}[c]{@{}c@{}}
\includegraphics[scale=.95]{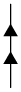}
\end{tabular}
= 2
\begin{tabular}[c]{@{} c@{}}
\includegraphics[scale=.95]{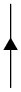}
\end{tabular}$}
\qquad
\subcaptionbox{\label{Fig098--Fig099}}{$\begin{tabular}[c]{@{}c@{}}
\includegraphics[scale=.95]{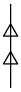}
\end{tabular}
= 2
\begin{tabular}[c]{@{} c@{}}
\includegraphics[scale=.95]{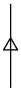}
\end{tabular}$}
\caption{Additional relations of $\widehat{\P}_{k}^{LR}(\Omega)$.}\label{Fig096--Fig099}
\end{figure}

\begin{example}
Figure~\ref{Fig101--Fig102} depicts multiplication of three diagrams in $\widehat{\P}_{6}^{LR}(\Omega)$ and Figure~\ref{Fig103--Fig104} shows an example where each of the diagrams and their product have $\a$-value 1.  Again, we use the dotted line to emphasize that the two closed decorations on the leftmost propagating edge belong to distinct blocks.
\end{example}

\begin{figure}[!ht]
\centering
\begin{tabular}[c]{l}
\includegraphics[scale=.95]{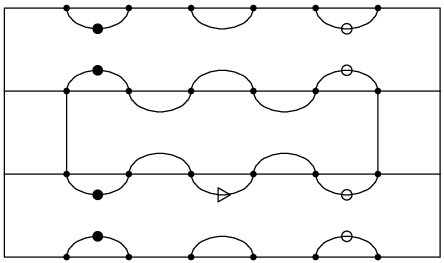}
\end{tabular} = 2 \begin{tabular}[c]{l}
\includegraphics[scale=.95]{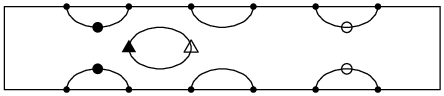}
\end{tabular}
\caption{Example of multiplication in $\widehat{\P}_{k}^{LR}(\Omega)$.}\label{Fig101--Fig102}
\end{figure} 

\begin{figure}[!ht]
\centering
\begin{tabular}[c]{l}
\includegraphics[scale=.95]{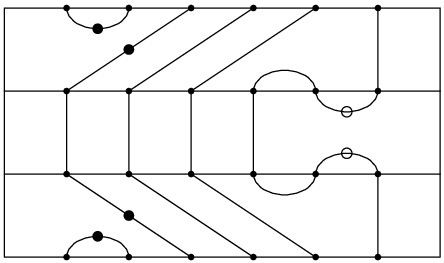}
\end{tabular}
= \begin{tabular}[c]{l}
\includegraphics[scale=.95]{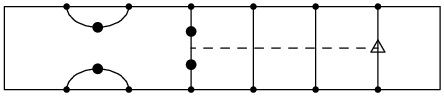}
\end{tabular}
\caption{Example of multiplication in $\widehat{\P}_{k}^{LR}(\Omega)$ with diagrams having $\a$-value 1.}\label{Fig103--Fig104}
\end{figure}

\end{subsection}


\begin{subsection}{Irreducible LR-decorated diagrams as a basis}

We need to show that a basis for $\widehat{\P}_{k}^{LR}(\Omega)$ consists of the set of LR-decorated diagrams having maximal blocks corresponding to nonidentity basis elements in $\V$.  That is, no block may contain adjacent decorations of the same type (open or closed).  To accomplish this task, we will make use of a diagram algebra version of Bergman's Diamond Lemma~\cite{Bergman.G:A}.  For other examples of this type of application of Bergman's Diamond Lemma, see~\cite{Green.R;Martin.P;Parker.A:A} and~\cite{Martin.P;Green.R;Parker.A:A}.

Define the function $r: T_{k}^{LR}(\Omega) \to T_{k}(\emptyset)$ via
\[
r(d)=d \text{ with all decorations and loops removed}.
\]
In the literature, if $d$ has no loops, then $r(d)$ is sometimes referred to as the ``shape'' of $d$.

Next, define a function $h: T_{k}^{LR}(\Omega) \to \Z^{+}\cup \{0\}$ via
\[
h(d)=\text{sum of the number of decorations and the number of loops}.
\]
Define $\leq_{\widehat{\P}}$ on $T_{k}^{LR}(\Omega)$ via $d < _{\widehat{\P}} d'$ if and only if $r(d)=r(d')$ and $h(d) < h(d')$.

Consider the collection of reductions determined by the relations of $\widehat{\P}_{k}^{LR}(\Omega)$ given in Definition~\ref{def:big diagram alg}.  If we apply any single reduction (loop removal or any other local reduction) to a diagram from $\widehat{\P}_{k}^{LR}(\Omega)$, then we obtain a scalar multiple of a strictly smaller diagram with respect to $\leq_{\widehat{\P}}$.  Thus, our reduction system (i.e., diagram relations) is compatible with $\leq_{\widehat{\P}}$.  Now, suppose that $d <_{\widehat{\P}} d'$ and let $d''$ be any other element from $\widehat{\P}_{k}^{LR}(\Omega)$.  Then $r(d''d)=r(d''d')$ and $r(dd'')=r(d'd'')$.  Since $r(d)=r(d')$, multiplying $d$ or $d'$ on the same side by $d''$ will increase the number of decorations and number of loops by the same amount.  So, we have $h(dd'') < h(d'd'')$ and $h(d''d) < h(d''d')$.  Therefore, $dd'' <_{\widehat{\P}} d'd''$ and $d''d <_{\widehat{\P}} d''d'$.  This shows that $\leq_{\widehat{\P}}$ is a semigroup partial order on $T_{k}^{LR}(\Omega)$.  Clearly, $\leq_{\widehat{\P}}$ satisfies the descending chain condition.

\begin{proposition}\label{prop:diagram independence}
The set of LR-decorated diagrams having no relations to apply forms a basis for $\widehat{\P}_{k}^{LR}(\Omega)$.
\end{proposition}

\begin{proof}
Let $\leq_{\widehat{\P}}$ be as above.  Following the setup of Bergman's Diamond Lemma, it remains to show that all of the ambiguities are resolvable.  

By inspecting the relations of Definition~\ref{def:big diagram alg}, we see that there are no inclusion ambiguities, so we only need to check that the overlap ambiguities are resolvable.  

Let $d$ be a diagram from $\widehat{\P}_{k}^{LR}(\Omega)$ and suppose that there are two competing reductions that we could apply.  If both reductions involve the same non-loop edge, then the ambiguity is easily seen to be resolvable since the algebra $\V$ is associative.  In particular, in the $\a$-value $1$ case, the reductions could involve two distinct blocks on the same edge, in which case, the order that we apply the reductions is immaterial.  If the reductions involve distinct edges, loop edges or otherwise, the ambiguity is quickly seen to be resolvable since the reductions commute.  Finally, suppose that the two competing reductions involve the same loop edge.  There are three possibilities for this loop edge: (a) the loop is undecorated, (b) the loop carries only one type of decoration (open or closed), and (c) the loop carries both types of symbols.  Note that (a) cannot happen since then there could not have been two competing reductions involving this edge to apply.  If (b) occurs, then any ambiguity involving this loop edge (including removing the loop) is resolvable since multiplication of closed (respectively, open) decorations is commutative and associative.  Finally, assume (c) occurs.  Note that the nature of our relations prevents the complete elimination of closed (respectively, open) decorations from this loop edge.  Since all loop relations involve either undecorated loops or loops decorated with a single type of decoration, this loop edge can never be removed.  Since $\V$ is associative and none of the relations involve both decoration types at the same time, the ambiguity is easily seen to be resolvable since the reductions commute.  

According to Bergman's Diamond Lemma~\cite{Bergman.G:A}, we can conclude that the set of LR-decorated diagrams having no relations to apply is a basis, as desired.
\end{proof}

\end{subsection}

\end{section}


\begin{section}{The simple and admissible diagrams}\label{sec:simple and admissible}

In this section, we define the diagram algebra $\D_n$ as a subalgebra of $\widehat{\P}_{n+2}^{LR}(\Omega)$ that turns out to be a faithful diagrammatic representation of $\TL(\C_n)$.  We will be able to quickly conclude that there is a surjective homomorphism from $\TL(\C_n)$ to $\D_n$.  In the sequel to this paper~\cite{Ernst.D:D}, we show that this homomorphism is injective, thus showing that the algebras are isomorphic.  In the next section of this paper, we define the admissible diagrams and show that they are a basis for $\D_n$.  In fact, we will show that the image of each monomial basis element of $\TL(\C_n)$ is admissible.


\begin{subsection}{Simple diagrams}

Define the \emph{simple diagrams} $d_{1}, d_{2}, \dots, d_{n+1}$ as in Figure~\ref{Fig107--Fig109}.  Note that the simple diagrams are elements of the basis for $\widehat{\P}_{n+2}^{LR}(\Omega)$ described in Proposition~\ref{prop:diagram independence}.

\begin{figure}[!ht]
\centering
\begin{eqnarray*}
d_{1} & = &
\begin{tabular}[c]{@{} c@{}}
\includegraphics[scale=.95]{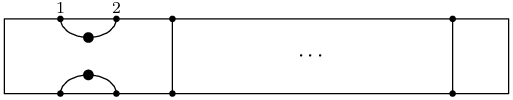}
\end{tabular}\\ 
& \vdots & \\
d_{i} & = &
\begin{tabular}[c]{@{} c@{}}
\includegraphics[scale=.95]{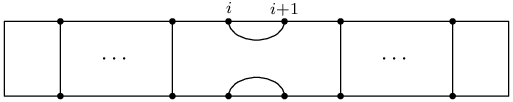}
\end{tabular}\\
& \vdots & \\
d_{n+1} & = &
\begin{tabular}[c]{@{} c@{}}
\includegraphics[scale=.95]{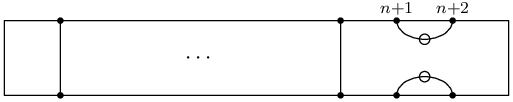}
\end{tabular}
\end{eqnarray*}
\caption{The simple diagrams.}\label{Fig107--Fig109}
\end{figure}

It is not immediately obvious, but we shall see that the algebra generated by the simple diagrams is infinite dimensional yet strictly smaller that $\widehat{\P}_{n+2}^{LR}(\Omega)$.

\begin{remark}\label{rem:affine C relations hold}
Checking that each of the following relations is satisfied for the simple diagrams is easily verified.
\begin{enumerate}
\item $d_{i}^{2}=\delta d_{i}$ for all $i$;
\item $d_{i}d_{j}=d_{i}d_{j}$ if $|i-j|>1$;
\item $d_{i}d_{j}d_{i}=d_{i}$ if $|i-j|=1$ and $1< i,j < n+1$;
\item $d_{i}d_{j}d_{i}d_{j}=2d_{i}d_{j}$ if $\{i,j\}=\{1,2\}$ or $\{n,n+1\}$.
\end{enumerate}
This shows that the simple diagrams satisfy the relations of $\TL(\C_n)$ given in Definition~\ref{def:TL(C)}.
\end{remark}

Finally, we are ready to define the diagram algebra that we are ultimately interested in.  Defining the algebra is easy, but having a description of a collection of basis diagrams is not.  The issue of the basis will be handled in Section~\ref{sec:basis}.

\begin{definition}\label{def:D_n}
Let $\D_{n}$ be the $\Z[\delta]$-subalgebra of $\widehat{\P}_{n+2}^{LR}(\Omega)$ generated as a unital algebra by $d_{1}, d_{2}, \dots, d_{n+1}$ with multiplication inherited from $\widehat{\P}_{n+2}^{LR}(\Omega)$.
\end{definition}

Now, define $\theta: \TL(\C_{n}) \to \D_{n}$ to be the function determined by $\theta(b_{i})=d_{i}$.  The next theorem follows quickly.

\begin{proposition}\label{prop:surjective homomorphism}
The map $\theta$ defined above is a surjective algebra homomorphism.
\end{proposition}

\begin{proof}
By Remark~\ref{rem:affine C relations hold}, the simple diagrams satisfy the relations of $\TL(\C_n)$.  This shows that $\theta$ is an algebra homomorphism, but since the simple diagrams generate $\D_n$, $\theta$ is surjective.
\end{proof}

The main result of the sequel to this paper~\cite{Ernst.D:D} is that $\theta$ is injective.

\end{subsection}


\begin{subsection}{Admissible diagrams}\label{subsec:admissible}

The next definition describes the set of $\C$-admissible diagrams, which will turn out to form a basis for $\D_n$.  Our definition of $\C$-admissible is motivated by the definition of $B$-admissible (after an appropriate change of basis) given by R.M.~Green in~\cite[Definition~2.2.4]{Green.R:H} for diagrams in the context of type $B$.  Since the Coxeter graph of type $\C$ is type $B$ at ``both ends'', the general idea is to build the axioms of $B$-admissible into our definition of $\C$-admissible on the left and right sides of our diagrams. 

\begin{definition}\label{def:admissible}
Let $d$ be an irreducible LR-decorated diagram.  Then we say that $d$ is \emph{$\C$-admissible}, or simply \emph{admissible}, if the following axioms are satisfied.
\begin{enumerate}[label=\rm{(C\arabic*)}]

\item \label{C1}The only loops that may appear are equivalent to the one in Figure~\ref{Fig110}.

\begin{figure}[!ht]
\centering
\includegraphics[scale=.95]{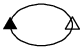}
\caption{The only allowable loop in $\C$-admissible diagrams.}\label{Fig110}
\end{figure}

\item \label{C2} Assume $\a(d)>1$ and let $e$ be the edge connected to node 1.  If $e$ is not connected to node $1'$, then it is decorated and the first decoration is a $\bcirc$.  If $e$ is connected to $1'$, then exactly one of the following three conditions are met:
\begin{enumerate}
\item $e$ is undecorated.
\item $e$ is decorated by a single $\btri$.
\item $e$ is decorated by a single block of decorations consisting of an alternating sequence of closed and open decorations such that the first decoration is a $\bcirc$.
\end{enumerate}
We have analogous restrictions for nodes $1'$, $n+2$, and $(n+2)'$, where we replace first with last for nodes $1'$ and $(n+2)'$ and closed decorations are replaced with open decorations for nodes $n+2$ and $(n+2)'$.

\item \label{C3} Assume $\a(d)=1$.  Then the western end of $d$ is equal to one of the diagrams in Figure~\ref{Fig110--Fig115}, where $u\in\{\emptyset,\btri\}$ and the other rectangles represent a sequence of blocks (possibly empty) such that each block is a single $\btri$.  Moreover, if $d$ is the diagram in Figure~\ref{Fig112}, then no more decorations occur on $d$.  Also, the occurrences of the $\bcirc$ decorations occurring on the propagating edges have the highest (respectively, lowest) relative vertical position of all decorations occurring on any propagating edge.  We have an analogous restrictions for the eastern end of $d$, where the closed decorations are replaced with open decorations.

\item \label{C4} No other $\bcirc$ or $\wcirc$ decorations appear on $d$ other than those required in~\ref{C2} and~\ref{C3}.

\end{enumerate}
Let $\Diag^{b}_{n}(\Omega)$ denote the set of all $\C$-admissible $(n+2)$-diagrams.
\end{definition}

\begin{figure}[!ht]
\centering
\subcaptionbox{\label{Fig111}}{\includegraphics[scale=.95]{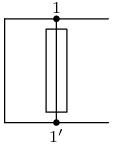}}
\qquad
\subcaptionbox{\label{Fig112}}{\includegraphics[scale=.95]{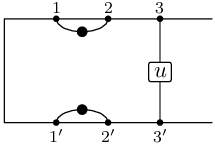}}
\qquad
\subcaptionbox{\label{Fig113}}{\includegraphics[scale=.95]{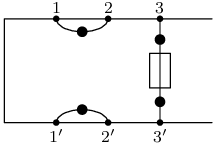}}\\
\subcaptionbox{\label{Fig114}}{\includegraphics[scale=.95]{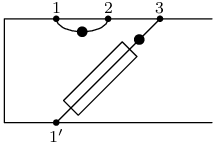}}
\qquad
\subcaptionbox{\label{Fig115}}{\includegraphics[scale=.95]{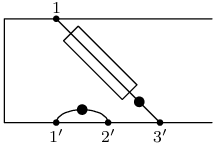}}
\caption{The western end of diagrams exhibiting axiom~\ref{C3}.}\label{Fig110--Fig115}
\end{figure}

\begin{remark}\label{rem:admissible}
We collect several comments concerning the admissible diagrams.
\begin{enumerate}
\item The only time an admissible diagram $d$ can have an edge adorned with both open and closed decorations is if $d$ is undammed (which only happens when $n$ is even) or if $d$ has a single propagating edge (which only happens when $n$ is odd).  The latter case coincides with part (c) of axiom~\ref{C2}. See parts (a) and (c) of Example~\ref{ex:decorated diagrams} for examples of diagrams having edges adorned with both types of decorations.

\item \label{rem:alternating decorations a=1} If $d$ is an admissible diagram with $\a(d)=1$, then the restrictions on the relative vertical position of decorations on propagating edges along with axiom~\ref{C3} imply that the relative vertical positions of closed decorations on the leftmost propagating edge and open decorations on the rightmost propagating edge must alternate.  In particular, the number of closed decorations occurring on the leftmost propagating edge differs from the number of open decorations occurring on the rightmost propagating edge by at most 1.  For example, if $d$ is the diagram in Figure~\ref{Fig116}, where the leftmost propagating edge carries $k$ $\btri$ decorations, then the rightmost propagating edge must carry $k$ $\wtri$ decorations, as well. 

\begin{figure}[!ht]
\centering
\includegraphics[scale=.95]{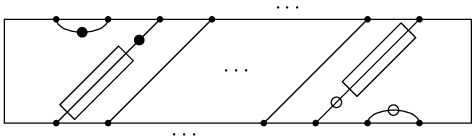}
\caption{Example of a diagram exhibiting axiom~\ref{C3}.}\label{Fig116}
\end{figure}

\item It is clear that $\Diag^{b}_{n}(\Omega)$ is an infinite set.  If an admissible diagram $d$ is undammed, then there is no limit to the number of loops given in axiom~\ref{C1} that may occur.  Also, if $d$ is an admissible diagram with exactly one propagating edge, then there is no limit to the width of the block of decorations that may occur on the lone propagating edge.  Furthermore, if $d$ is admissible with $\a(d)=1$, then there is no limit to the number of $\btri$-blocks (respectively, $\wtri$-blocks) that may occur on the leftmost (respectively, rightmost) propagating edge.  

\item \label{rem:admissible basis} Each of the admissible diagrams is a basis element of $\widehat{\P}_{n+2}^{LR}(\Omega)$.

\item The symbol $b$ in the notation $\Diag^{b}_{n}(\Omega)$ is to emphasize that we are constructing a set of diagrams that is intended to correspond to the monomial basis of $\TL(\C_{n})$.  In a sequel to this paper, we will construct diagrams that correspond to the ``canonical basis'' of $\TL(\C_{n})$, which is defined for arbitrary Coxeter groups in~\cite{Green.R;Losonczy.J:A}.
\end{enumerate}
\end{remark}

\begin{definition}\label{def:M[Diag_n^b(V)]}
Let $\mathcal{M}[\Diag_{n}^{b}(\Omega)]$ be the $\Z[\delta]$-submodule of $\widehat{\P}_{n+2}^{LR}(\Omega)$ spanned by the admissible diagrams.
\end{definition}

\begin{proposition}\label{prop:admissible basis}
The set of admissible diagrams $\Diag_{n}^{b}(\Omega)$ is a basis for the module $\mathcal{M}[\Diag_{n}^{b}(\Omega)]$.
\end{proposition}

\begin{proof}
Linear independence follows immediately from Remark~\ref{rem:admissible}(\ref{rem:admissible basis}).
\end{proof}

\end{subsection}


\begin{subsection}{Temperley--Lieb diagram algebras of type $B$}\label{subsec:type B}

We will briefly discuss how $\TL(B_{n})$ and $\TL(B'_{n})$ are related to $\D_{n}$.

\begin{definition}\label{def:DTL(B_n)}
Let $\DTL(B_{n})$ and $\DTL(B'_{n})$ denote the subalgebras of $\D_{n}$ generated by the simple diagrams $d_{1}, d_{2}, \dots, d_{n}$ and $d_{2}, d_{3}, \dots, d_{n+1}$, respectively.  We refer to $\DTL(B_{n})$ (respectively, $\DTL(B'_{n})$) as the \emph{Temperley--Lieb diagram algebra of type $B$} (respectively, \emph{type $B'$}). 
\end{definition}

It is clear that $\DTL(B_{n})$ (respectively, $\DTL(B'_{n})$) consists entirely of diagrams that are decorated with only closed (respectively, open) decorations.  Also, note that all of the technical requirements about how to decorate a diagram $d$ when $\a(d)=1$ are irrelevant since only the leftmost (respectively, rightmost) propagating edge can carry decorations in  $\DTL(B_{n})$ (respectively, $\DTL(B'_{n})$).  The following fact is implicit in~\cite[\textsection 2]{Green.R:H} after the appropriate change of basis involving a change of basis for the decoration set.

\begin{proposition}
As $\Z[\delta]$-algebras, $\TL(B_{n}) \cong \DTL(B_{n})$ and $\TL(B'_{n}) \cong \DTL(B'_{n})$, where each isomorphism is determined by $b_{i} \mapsto d_{i}$ for the appropriate restrictions on $i$.  \hfill $\qed$
\end{proposition}

After making the appropriate change of basis on the decoration set (which involves making a change of basis on the rank 3 Verlinde algebra), the basis diagrams in $\DTL(B_{n})$ (respectively, $\DTL(B'_{n})$) become $B$-admissible in the sense of~\cite{Green.R:H, Green.R:M}.  Moreover, it is easily verified that the axioms for $B$-admissible given in~\cite[Definition~2.2.4]{Green.R:H} imply (again, under the appropriate change of basis involving the decoration set) that all of the basis diagrams in $\DTL(B_{n})$ and $\DTL(B'_{n})$ are $\C$-admissible.  

\end{subsection}

\end{section}


\begin{section}{A basis for $\D_{n}$}\label{sec:basis}

Our main objective in the remainder of this paper is to show that the admissible diagrams form a basis for $\D_n$.  Before proceeding, we wish to outline our method of attack.  We will show the following:

\begin{enumerate}
\item \label{hard1} The admissible diagrams are generated by the simple diagrams (Proposition~\ref{prop:admissibles gen by simples}).
\item \label{hard2} The module $\mathcal{M}[\Diag_{n}^{b}(\Omega)]$ is closed under multiplication, making it a subalgebra of $\widehat{\P}^{LR}_{n+2}(\Omega)$ (Corollary~\ref{cor:module is subalgebra}).
\item \label{easy} The algebras $\mathcal{M}[\Diag_{n}^{b}(\Omega)]$ and $\D_n$ are equal having the admissible diagrams as a basis (Theorem~\ref{thm:module admissibles is D_n}).
\end{enumerate}

Items (\ref{hard1}) and (\ref{hard2}) above require numerous technical lemmas.  However, once we overcome these difficulties, (\ref{easy}) will yield itself easily.


\begin{subsection}{Preparatory lemmas}
Our first significant obstacle in proving that the $\C$-admissible diagrams form a basis for $\D_{n}$ is proving that the admissible diagrams are generated by the simple diagrams (see Proposition~\ref{prop:admissibles gen by simples}).  To achieve this end, we require several intermediate results.

If $d$ is an admissible diagram, then we say that a non-propagat\-ing edge joining $i$ to $i+1$ (respectively, $i'$ to $(i+1)'$) is \emph{simple} if it is identical to the edge joining the same nodes of the simple diagram $d_i$.  Note that a simple edge is undecorated except when one of the vertices is 1 or $1'$ (respectively, $n+2$ or $(n+2)'$), in which case it is decorated by only a single $\bcirc$ (respectively, $\wcirc$).  

The next six lemmas mimic Lemmas~5.1.4--5.1.7 in~\cite{Green.R:H}.  The proof of each lemma is immediate and throughout we assume that $d$ is admissible.  Figure~\ref{Fig117--Fig122} provides visual representations of each lemma, where $x$ represents an arbitrary (possibly empty) block of decorations.  Each of Lemmas~\ref{lem:nonprops connect adj verts}--\ref{lem:move decoration 2} have left-right symmetric analogues (perhaps involving closed decorations), as well as versions that involve edges in the south face.

\begin{figure}[!ht]
\centering
\subcaptionbox{Lemma~\ref{lem:nonprops connect adj verts}\label{Fig117}}{\includegraphics[scale=.95]{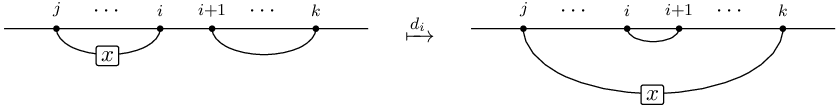}}\\
\vspace{2ex}
\subcaptionbox{Lemma~\ref{lem:nonprops connect adj verts 2}\label{Fig118}}{\includegraphics[scale=.95]{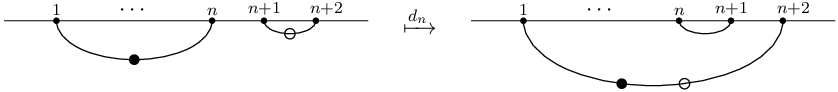}}\\
\vspace{2ex}
\subcaptionbox{Lemma~\ref{lem:swap prop nonprop}\label{Fig119}}{\includegraphics[scale=.95]{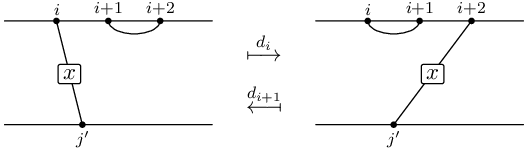}}\\
\vspace{2ex}
\subcaptionbox{Lemma~\ref{lem:make decoration}\label{Fig120}}{\includegraphics[scale=.95]{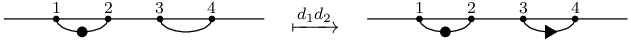}}\\
\vspace{2ex}
\subcaptionbox{Lemma~\ref{lem:move decoration}\label{Fig121}}{\includegraphics[scale=.95]{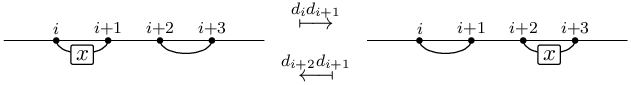}}\\
\vspace{2ex}
\subcaptionbox{Lemma~\ref{lem:move decoration 2}\label{Fig122}}{\includegraphics[scale=.95]{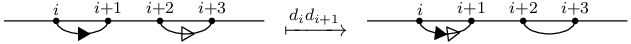}}\\
\vspace{2ex}
\caption{Visual representations of Lemmas~\ref{lem:nonprops connect adj verts}--\ref{lem:move decoration 2}}\label{Fig117--Fig122}
\end{figure}

\begin{lemma}\label{lem:nonprops connect adj verts}
Assume that in the north face of $d$ there is an edge, say $e$, connecting node $j$ to node $i$, and assume that there is another, undecorated, edge, say $e'$, connecting node $i+1$ to node $k$ with $j<i$ and $i+1<k<n+2$.    Then $d_{i}d$ is the admissible diagram that results from $d$ by removing $e'$, disconnecting $e$ from node $i$ and reattaching it to node $k$, and adding a simple edge to $i$ and $i+1$ (note that edge $e$ maintains its original decorations). See Figure~\ref{Fig117}. \hfill $\qed$
\end{lemma}

\begin{lemma}\label{lem:nonprops connect adj verts 2}
Assume that in the north face of $d$ there is an edge, say $e$, connecting node $1$ to node $n$ labeled by a single $\bcirc$ (this can happen only if $n$ is even), and assume that there is a simple edge, say $e'$, connecting node $n+1$ to node $n+2$ (which must be labeled by a single $\wcirc$).  Then $d_{n}d$ is the admissible diagram that results from $d$ by joining the right end of $e$ to the left end of $e'$, and adding a simple edge that joins $n$ to $n+1$.  Note that the new edge formed by joining $e$ and $e'$ connects node $1$ to node $n+2$ and is labeled by the block $\bcirc \wcirc$.  See Figure~\ref{Fig118}. \hfill $\qed$
\end{lemma}

\begin{lemma}\label{lem:swap prop nonprop}
Assume that $d$ has a propagating edge, say $e$, joining node $i$ to node $j'$ with $1<i<n$.  Further, assume that there is a simple edge, say $e'$, joining nodes $i+1$ and $i+2$.  Then $d_{i}d$ is the admissible diagram that results from $d$ by removing $e'$, disconnecting $e$ from node $i$ and reattaching it to node $i+2$, and adding a simple edge to $i+1$ and $i+2$ (note that $e$ retains its original decorations).  See Figure~\ref{Fig119}. \hfill $\qed$
\end{lemma}

\begin{lemma}\label{lem:make decoration}
Assume that $d$ has simple edges joining node $1$ to node $2$ and node $3$ to node $4$.  Then $d_{1}d_{2}d$ is the admissible diagram that results from $d$ by adding a $\btri$ to the edge joining $3$ to $4$.  See Figure~\ref{Fig120}. \hfill $\qed$
\end{lemma}

\begin{lemma}\label{lem:move decoration}
Assume that $d$ has two edges, say $e$ and $e'$, joining node $i$ to node $i+1$ and node $i+2$ to node $i+3$, respectively, where $1<i<n-1$ and $e'$ is simple.  Then $d_{i}d_{i+1}d$ is the admissible diagram that results from $d$ by removing the decorations from $e$ and adding them to $e'$.  This procedure has an inverse, since $d_{i+2}d_{i+1}(d_{i}d_{i+1}d)=d$.  See Figure~\ref{Fig121}.  \hfill $\qed$
\end{lemma}

\begin{lemma}\label{lem:move decoration 2}
Assume that $d$ has two edges, say $e$ and $e'$, joining node $i$ to node $i+1$ and node $i+2$ to node $i+3$, respectively, with $1<i<n-1$. Further, assume that $e$ is decorated by a single $\btri$ decoration only and that $e'$ is decorated by a single $\wtri$ decoration only.  Then $d_{i+2}d_{i+1}d$ is the admissible diagram that results from $d$ by removing the $\wtri$ decoration from $e'$ and adding it to $e$ to the right of the $\btri$ decoration. See Figure~\ref{Fig122}.  \hfill $\qed$
\end{lemma}

\end{subsection}


\begin{subsection}{The admissible diagrams are generated by the simple diagrams}

Next, we state and prove several lemmas that we will use to prove that each admissible diagram can be written as a product of simple diagrams in $\D_{n}$.

\begin{lemma}\label{lem:a-value=1 diagrams gen by simples}
If $d$ is an admissible diagram with $\a(d)=1$, then $d$ can be written as a product of simple diagrams.
\end{lemma}

\begin{proof}
Assume that $d$ is an admissible diagram with $\a(d)=1$.  The proof is an exhaustive case by case check, where we consider all the possible diagrams that are consistent with axiom~\ref{C3}.  We consider five cases (see Figure~\ref{Fig163--Fig167}); all remaining cases follow by analogous arguments.  

\begin{figure}[!ht]
\centering
\subcaptionbox{\label{Fig163}}{\includegraphics[scale=.95]{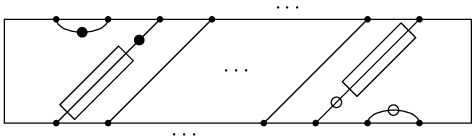}}
\qquad
\subcaptionbox{\label{Fig164}}{\includegraphics[scale=.95]{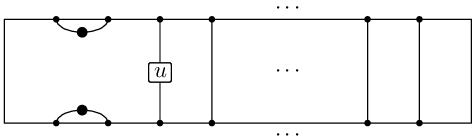}}\\
\subcaptionbox{\label{Fig165}}{\includegraphics[scale=.95]{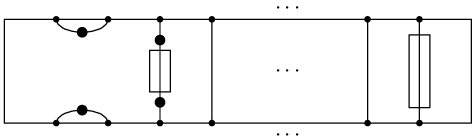}}
\qquad
\subcaptionbox{\label{Fig166}}{\includegraphics[scale=.95]{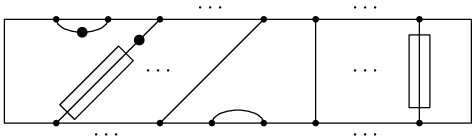}}
\subcaptionbox{\label{Fig167}}{\includegraphics[scale=.95]{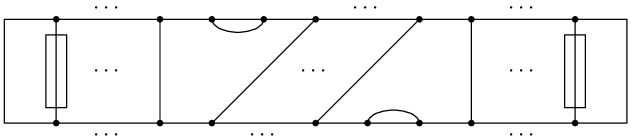}}
\caption{The five cases of Lemmas~\ref{lem:a-value=1 diagrams gen by simples} and~\ref{lem:closed under mult 7}.}\label{Fig163--Fig167}
\end{figure}

Case (1). First, assume that $d$ is the diagram in Figure~\ref{Fig163}, where the leftmost propagating edge carries $k$ $\btri$ decorations, and hence, the rightmost propagating edge carries $k$ $\wtri$ decorations by  Remark~\ref{rem:admissible}(\ref{rem:alternating decorations a=1}).  In this case, it can quickly be verified that we can obtain $d$ via
\[
d=(d_{z_{1}}d_{z_{2}})^{k}d_{z_{1}}d_{n+1},
\]
where $d_{z_{1}}=d_{1}d_{2}\cdots d_{n}$ and $d_{z_{2}}=d_{n+1}d_{n}\cdots d_{2}$.  Therefore, $d$ can be written as a product of simple diagrams, as desired.

Case (2). For the second case, assume that $d$ is the diagram in Figure~\ref{Fig164}, where $u\in\{\emptyset,\btri\}$.  Note that $d$ does not carry any open decorations.  In this case, either $d=d_1$ or $d=d_{1}d_{2}d_{1}$, and regardless $d$ can be written as a product of simple diagrams, as expected.

Case (3). For the third case, assume that $d$ is the diagram in Figure~\ref{Fig165}, where the leftmost propagating edge carries $k-1$ $\btri$ decorations, so that the rightmost propagating edge carries $k$ $\wtri$ decorations.  Then
\[
d=(d_{z_{1}}d_{z_{2}})^{k}d_{1},
\]
where $d_{z_{1}}$ and $d_{z_{2}}$ are as in case (1), and hence $d$ can be written as a product of simple diagrams.

Case (4). Next, assume that $d$ is the diagram in Figure~\ref{Fig166}, where the simple edge in the south face connects nodes $j'$ and $(j+1)'$ with $1'<j'<(n+1)'$, and the leftmost propagating edge carries $k$ $\btri$ decorations.  Then by Remark~\ref{rem:admissible}(\ref{rem:alternating decorations a=1}), the rightmost propagating edge carries $l$ $\wtri$ decorations, where $l=k$ or $k+1$.  If $l=k+1$, then define $d'$ to be the diagram in Figure~\ref{Fig163}, where the leftmost (respectively, rightmost) propagating edge carries $k$ $\btri$ (respectively, $\wtri$) decorations.  By case (1), $d'$ can be written as a product of simple diagrams.  We see that
\[	
d=d' d_{n}d_{n-1}\cdots d_{j+1}d_{j},
\]
which implies that $d$ can be written as a product of simple diagrams, as desired.  If, on the other hand, $l=k$, then define $d'$ to be identical to $d$ except that the last $\btri$ decoration occurring on the leftmost propagating edge has been removed.  Then by the subcase we just completed (where the rightmost propagating edge carried one more $\wtri$ decoration than the leftmost propagating edge carried $\btri$ decorations), $d'$ can be written as a product of simple diagrams.  We see that
\[
d=d' d_{j-1}d_{j-2}\cdots d_{2}d_{1}d_{2} \cdots d_{j-1}d_{j},
\]
which implies that $d$ can be written as a product of simple diagrams.

Case (5).  For the final case, assume that $d$ is the diagram in Figure~\ref{Fig167}, where the simple edge in the north face connects nodes $i$ and $i+1$ with $1<i<n+1$, the simple edge in the south face connects nodes $j'$ and $(j+1)'$ with $1'<j'<(n+1)'$, and the leftmost propagating edge carries $k$ $\btri$ decorations.  Then again by Remark~\ref{rem:admissible}(\ref{rem:alternating decorations a=1}), the rightmost propagating edge carries $l$ $\wtri$ decorations, where $|k-l| \leq 1$.  Without loss of generality, assume that $k \leq l$, so that $l=k$ or $k+1$.  If $l=k+1$, define $d'$ to be the diagram in Figure~\ref{Fig166}, where the leftmost propagating edge carries $k-1$ $\btri$ decorations while the rightmost propagating edge carries $k$ $\wtri$ decorations.  By case (4), $d'$ can be written as a product of simple diagrams.  We see that
\[
d=d_{i}d_{i+1}\cdots d_{n}d_{n+1}d_{n}\cdots d_{3}d_{2}d',
\]
which implies that $d$ can be written as a product of simple diagrams.  If, on the other hand, $k=l$, then without loss of generality, assume that the first decoration occurring on the leftmost propagating edge has the highest relative vertical position of all decorations occurring on propagating edges.  Define $d'$ to be the diagram in Figure~\ref{Fig116}, where the leftmost (respectively, rightmost) propagating edge carries $k-1$ $\btri$ (respectively, $\wtri$) decorations.  Again, by case (4), $d'$ can be written as a product of simple diagrams.  Also, we see that
\[
d=d_{i}d_{i-1}\cdots d_{3}d_{2} d' d_{n}d_{n-1} \cdots d_{j+1}d_{j},
\]
which implies that $d$ can be written as a product of simple diagrams, as desired.
\end{proof}

\begin{lemma}\label{lem:a-value=max dammed diagrams gen by simples}
If $d$ is an admissible diagram with $1< \a(d)< \lfloor \frac{n+2}{2} \rfloor$ such that all non-propagating edges are simple, then $d$ can be written as a product of simple diagrams.
\end{lemma}

\begin{proof}
Let $d$ be an admissible diagram with $1< \a(d)< \lfloor \frac{n+2}{2} \rfloor$ such that all non-propagating edges are simple.  (Note that the restrictions on $\a(d)$ imply that $d$ has more than one propagating edge and at least one non-propagating edge.)   We consider two cases, where the second case has two subcases.

Case (1). First, assume that $d$ has a vertical propagating edge, say $e_{i}$, joining $i$ to $i'$.  Now, define the admissible diagrams $d'$ and $d''$ to be the diagrams in Figures~\ref{Fig131} and~\ref{Fig132}, respectively, where each of the shaded regions is identical to the corresponding regions of $d$.  Then $d=d'd''$.  Furthermore, since $d'$ (respectively, $d''$) satisfies requirement~\ref{D1} for LR-decorated diagrams, the diagram is only decorated with closed (respectively, open) decorations.  Since $d$ is admissible, $d' \in \DTL(B_{n})$ while $d'' \in \DTL(B'_{n})$.  This implies that both $d'$ and $d''$ can be written as a product of simple diagrams.  Therefore, $d$ can be written as a product of simple diagrams, as desired.

\begin{figure}[!ht]
\centering
\subcaptionbox{\label{Fig131}}{\includegraphics[scale=.95]{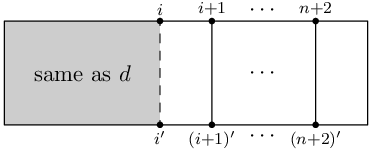}}
\qquad
\subcaptionbox{\label{Fig132}}{\includegraphics[scale=.95]{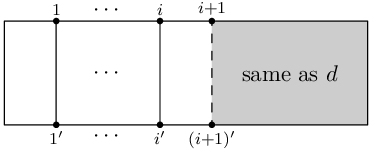}}
\caption{Diagrams for case (1) of the proof of Lemma~\ref{lem:a-value=1 diagrams gen by simples}.}\label{Fig131--Fig132}
\end{figure}

Case (2).  Next, assume that $d$ has no vertical propagating edges.  Suppose that the leftmost propagating edge joins node $i$ in the north face to node $j'$ in the south face, and without loss of generality, assume that $j<i$.  (Note that since $d$ has more than one propagating edge, $i<n+2$.)  We wish to make use of case (1), but we must consider two subcases.

(a) For the first subcase, assume that $j \neq 1$.  Since $d$ is admissible, $d$ must be the diagram in Figure~\ref{Fig133},where $x$ on the propagating edge from $i$ to $j'$ is either trivial (i.e., the edge is undecorated) or equal to a single $\btri$ decoration.  Define the admissible diagram $d'$ to be the diagram in Figure~\ref{Fig134}, where the leftmost propagating edge carries the same decoration as the leftmost propagating edge in $d$ and the shaded region is identical to the corresponding region of $d$.  By case (1), $d' $ can be written as a product of simple diagrams.  By making repeated applications of Lemma~\ref{lem:swap prop nonprop}, we can transform $d'$ into $d$, which shows that $d$ can be written as a product of simple diagrams, as desired.

\begin{figure}[!ht]
\centering
\begin{subfigure}[c]{0.05\textwidth}
\caption{}\label{Fig133}
\end{subfigure}
\begin{minipage}[c]{.85\textwidth}
\includegraphics[scale=.95]{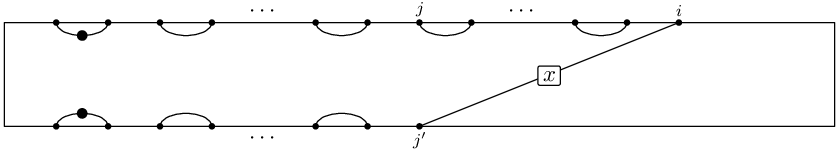}
\end{minipage}
\\[1ex]
\begin{subfigure}[c]{0.05\textwidth}
\caption{}\label{Fig134}
\end{subfigure}
\begin{minipage}[c]{.85\textwidth}
\includegraphics[scale=.95]{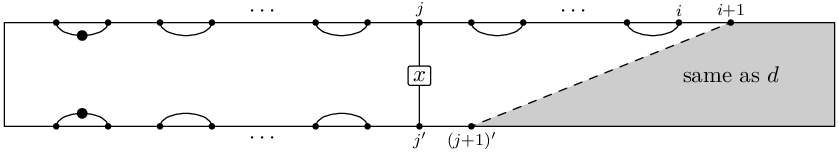}
\end{minipage}
\caption{Diagrams for case (2)(a) of the proof of Lemma~\ref{lem:a-value=1 diagrams gen by simples}.}\label{Fig133--Fig134}
\end{figure}

(b) For the second subcase, assume that $j=1$, so that $d$ is the diagram in Figure~\ref{Fig135}.  Since $1< \a(d)< \lfloor \frac{n+2}{2} \rfloor$, there is at least one other propagating edge occurring to the right of the leftmost propagating edge.  Furthermore, since the number of non-propagating edges in the north face is equal to the number of non-propagating edges in the south face, there is at least one undecorated non-propagating edge in the south face of $d$.  By making repeated applications, if necessary, of the southern version of Lemma~\ref{lem:swap prop nonprop}, we may assume that $d$ is the diagram in Figure~\ref{Fig136}.

\begin{figure}[!ht]
\centering
\begin{subfigure}[c]{0.05\textwidth}
\caption{}\label{Fig135}
\end{subfigure}
\begin{minipage}[c]{.58\textwidth}
\includegraphics[scale=.95]{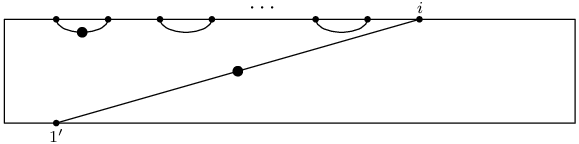}
\end{minipage}
\\[1ex]
\begin{subfigure}[c]{0.05\textwidth}
\caption{}\label{Fig136}
\end{subfigure}
\begin{minipage}[c]{.58\textwidth}
\includegraphics[scale=.95]{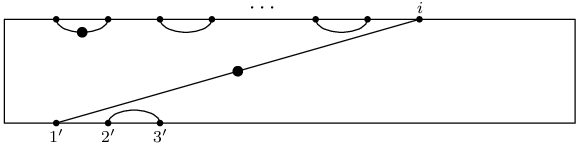}
\end{minipage}
\\[1ex]
\begin{subfigure}[c]{0.05\textwidth}
\caption{}\label{Fig137}
\end{subfigure}
\begin{minipage}[c]{.58\textwidth}
\includegraphics[scale=.95]{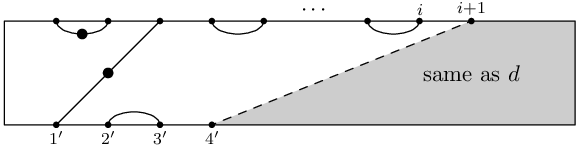}
\end{minipage}
\\[1ex]
\begin{subfigure}[c]{0.05\textwidth}
\caption{}\label{Fig138}
\end{subfigure}
\begin{minipage}[c]{.58\textwidth}
\includegraphics[scale=.95]{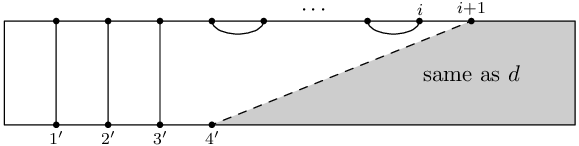}
\end{minipage}
\caption{Diagrams for case (2)(b) of the proof of Lemma~\ref{lem:a-value=1 diagrams gen by simples}.}\label{Fig135--Fig138}
\end{figure}

Now, define the admissible diagrams $d'$ and $d''$ via the diagrams in Figures~\ref{Fig137} and~\ref{Fig138}, respectively, where the shaded regions are identical to the corresponding regions of $d$.  By case (1), $d''$ can be written as a product of simple diagrams.  Also, we see that $d'=d_{1}d_{2}d''$, which implies that $d'$ can be written as a product of simple diagrams, as well.  By making repeated applications of Lemma~\ref{lem:swap prop nonprop}, we must have that $d$ can be written as a product of simple diagrams.
\end{proof}

\begin{lemma}\label{lem:a-value=max undammed diagrams gen by simples}
If $n$ is odd and $d$ is an admissible diagram with $\a(d)=\lfloor \frac{n+2}{2} \rfloor$ such that all non-propagating edges are simple, then $d$ can be written as a product of simple diagrams.
\end{lemma}

\begin{proof}
Assume that $n$ is odd and that $d$ is an admissible diagram with $\a(d)=\lfloor \frac{n+2}{2} \rfloor$.  In this case, $d$ has a unique propagating edge.  Also, assume that all of the non-propagating edges of $d$ are simple.  The proof is an exhaustive case by case check, where we consider the possible edges that are consistent with axioms~\ref{C2} and~\ref{C4} of Definition~\ref{def:admissible}.  We consider five cases (see Figure~\ref{Fig139--Fig145}); all remaining cases follow by analogous arguments.  

\begin{figure}[!ht]
\centering
\subcaptionbox{\label{Fig139}}{\includegraphics[scale=.95]{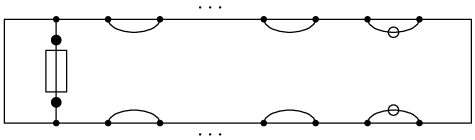}}
\qquad
\subcaptionbox{\label{Fig140}}{\includegraphics[scale=.95]{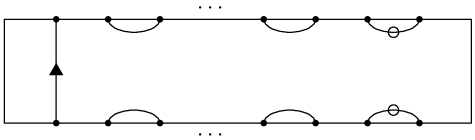}}
\subcaptionbox{\label{Fig141}}{\includegraphics[scale=.95]{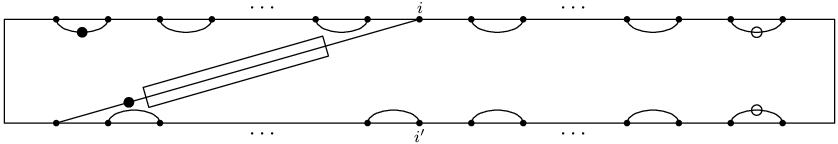}}
\subcaptionbox{\label{Fig143}}{\includegraphics[scale=.95]{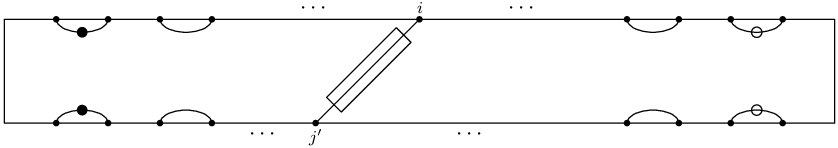}}
\subcaptionbox{\label{Fig145}}{\includegraphics[scale=.95]{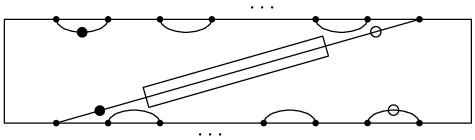}}
\caption{The five cases of Lemma~\ref{lem:a-value=max undammed diagrams gen by simples}.}\label{Fig139--Fig145}
\end{figure}

Case (1).  For the first case, assume that $d$ is the diagram in Figure~\ref{Fig139}, where the rectangle on the propagating edge is equal to a block consisting of an alternating sequence of $k-1$ $\btri$ decorations and $k$ $\wtri$ decorations.  It is quickly verified that
\[
d=(d_{\E}d_{\O})^{k}d_{\E},
\]
where
\[
d_{\E}=d_{2}d_{4}\cdots d_{n+1}
\]
and
\[
d_{\O}=d_{1}d_{3}\cdots d_{n}.
\]
This shows that $d$ can be written as a product of simple diagrams, as desired.

Case (2).  For the second case, assume that $d$ is the diagram in Figure~\ref{Fig140}.  In this case, we see that
\[
d=d_{2}d_{1}d_{2}d_{4}\cdots d_{n-1}d_{n+1},
\]
which shows that $d$ can be written as a product of simple diagrams.

Case (3).  Next, assume that $d$ is the diagram in Figure~\ref{Fig141}.  (Note that $i$ must be odd.)  If the rectangle is empty, then
\[
d=d_{1}d_{3}\cdots d_{i-2}d_{\E},
\]
where $d_{\E}$ is as in case (1).  On the other hand, if the rectangle is nonempty, so that the rectangle is equal to a block consisting of an alternating sequence of $k$ $\btri$ decorations and $l$ $\wtri$ decorations, where $l=k$ or $k+1$, define the admissible diagram $d'$ to be the one in Figure~\ref{Fig139}, where the rectangle on the propagating edge is equal to a block consisting of an alternating sequence of $k-1$ $\btri$ decorations and $k$ $\wtri$ decorations.  By case (1), $d'$ can be written as a product of simple diagrams.  If $k=l$, then we see that 
\[
d=d_{i-2}d_{i-4}\cdots d_{3}d_{1}d'.
\]
If, on the other hand, $l=k+1$, then we see that
\[
d=d_{i+1}d_{i+3}\cdots d_{n-1}d_{n+1}d_{\O}d',
\]
where $d_{\O}$ is as in case (1).  This shows that $d$ can be written as a product of simple diagrams.

Case (4).  Now, assume that $d$ is the diagram in Figure~\ref{Fig143}, where $i,j \notin \{1,n+2\}$ and the rectangle on the propagating edge is equal to a block consisting of an alternating sequence of $k$ $\btri$ decorations and $l$ $\wtri$ decorations with $|k-l| \leq 1$.  (Note that $i$ and $j$ must be odd.)  Without loss of generality, assume that $k \leq l$, so that $l=k$ or $k+1$.   Now, assume that the last decoration on the propagating edge is a $\btri$; the case with the last decoration being a $\wtri$ is handled with an analogous argument.  If $l=k$ (respectively, $l=k+1$), then the first decoration on the propagating edge is a $\wtri$ (respectively, $\btri$).  In either case, define the admissible diagram $d'$ via the diagram in Figure~\ref{Fig141}, where the rectangle on the propagating edge is equal to a block consisting of an alternating sequence of $k-1$ $\btri$ decorations and $l$ $\wtri$ decorations.  By case (3), $d'$ can be written as a product of simple diagrams.  Then it is quickly verified that
\[
d=d'd_{1}d_{3}\cdots d_{j-2}d_{j},
\]
and so $d$ can be written as a product of simple diagrams.

Case (5).  For the final case, assume that $d$ is the diagram of Figure~\ref{Fig145}, where the rectangle on the propagating edge is equal to a block consisting of an alternating sequence of $k$ $\btri$ decorations and $k$ $\wtri$ decorations.  It is quickly seen that
\[
d=(d_{\O}d_{\E})^{k+1},
\]
where $d_{\O}$ and $d_{\E}$ are as in case (1).  So, $d$ can be written as a product of simple diagrams, as expected.
\end{proof}

By stringing together the previous lemmas, we are able to prove the following proposition.

\begin{proposition}\label{prop:admissibles gen by simples}
Each admissible diagram can be written as a product of simple diagrams.  In particular, the admissible diagrams are contained in $\D_{n}$.
\end{proposition}

\begin{proof}
Let $d$ be an admissible diagram.  Lemma~\ref{lem:nonprops connect adj verts}, and if necessary Lemma~\ref{lem:nonprops connect adj verts 2}, along with their analogues, allow us to assume that all of the non-propagating edges of $d$ join adjacent vertices.  Furthermore, Lemmas~~\ref{lem:make decoration},~\ref{lem:move decoration}, and~\ref{lem:move decoration 2}, along with their analogues, allow us to assume that all of the non-propagating edges of $d$ are simple.  We now consider four distinct cases: (1) $\a(d)=1$, (2) $1< \a(d)< \lfloor \frac{n+2}{2} \rfloor$, (3) $\a(d)=\lfloor \frac{n+2}{2} \rfloor$ with $n$ odd (i.e., $d$ has a unique propagating edge), and (4) $\a(d)=\frac{n+2}{2}$ with $n$ even (i.e., $d$ is undammed).

\begin{figure}[!ht]
\centering
\begin{subfigure}[c]{0.05\textwidth}
\caption{}\label{Fig146}
\end{subfigure}
\begin{minipage}[c]{.68\textwidth}
\includegraphics[scale=.95]{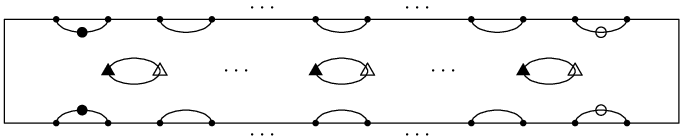}
\end{minipage}
\\[1ex]
\begin{subfigure}[c]{0.05\textwidth}
\caption{}\label{Fig147}
\end{subfigure}
\begin{minipage}[c]{.68\textwidth}
\includegraphics[scale=.95]{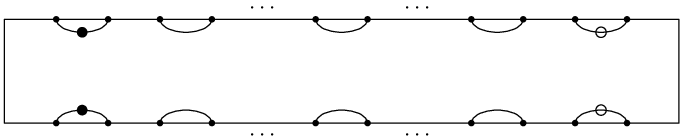}
\end{minipage}
\\[1ex]
\begin{subfigure}[c]{0.05\textwidth}
\caption{}\label{Fig148}
\end{subfigure}
\begin{minipage}[c]{.68\textwidth}
\includegraphics[scale=.95]{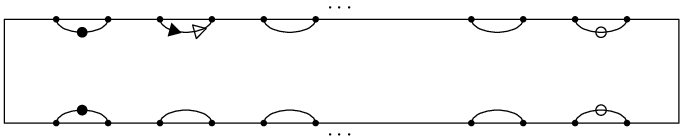}
\end{minipage}
\\[1ex]
\begin{subfigure}[c]{0.05\textwidth}
\caption{}\label{Fig149}
\end{subfigure}
\begin{minipage}[c]{.68\textwidth}
\includegraphics[scale=.95]{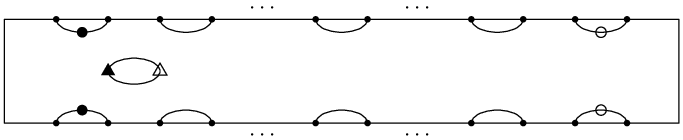}
\end{minipage}
\caption{Diagrams for case (4) of the proof of Proposition~\ref{prop:admissibles gen by simples}.}\label{Fig146--Fig149}
\end{figure}

Cases (1), (2), and (3) follow immediately from Lemmas~~\ref{lem:a-value=1 diagrams gen by simples},~\ref{lem:a-value=max dammed diagrams gen by simples}, and~\ref{lem:a-value=max undammed diagrams gen by simples}, respectively.  For the final case, assume that $\a(d)=\frac{n+2}{2}$ with $n$ even.  Then $d$ is undammed and based on our simplifying assumptions, we must have $d$ equal to the diagram in Figure~\ref{Fig146}, where there are $k$ loop edges (we allow $k=0$).  Define the admissible diagram 
\[
d_{\O}=d_{1}d_{3} \cdots d_{n+1}.
\]
Then $d_{\O}$ must be equal to the diagram in Figure~\ref{Fig147}.  In particular, $d_{\O}$ is identical to $d$, except that is has no loop edges.  If $d$ has no loop edges (i.e., $k=0$), then we are done.  Suppose $k>0$.  By making the appropriate repeated applications of the left and right-handed versions of Lemmas~~\ref{lem:make decoration} and~\ref{lem:move decoration} and a single application of Lemma~\ref{lem:move decoration 2}, there exists a sequence of simple diagrams $d_{i_{1}}, d_{i_{2}}, \dots, d_{i_{m}}$ such that $(d_{i_{1}}d_{i_{2}} \cdots d_{i_{m}})d_{\O}$ is equal to the diagram in Figure~\ref{Fig148}.  But then $d_{3}(d_{i_{1}}d_{i_{2}} \cdots d_{i_{m}})d_{\O}$ must be equal to the diagram in Figure~\ref{Fig149}.  To produce $k$ loops, we repeat this process $k-1$ more times.  That is,
\[
d=\left(d_{3}(d_{i_{1}}d_{i_{2}}\cdots d_{i_{m}})\right)^{k}d_{\O}.
\]
This shows that $d$ can be written as a product of simple diagrams, as desired.
\end{proof}

\end{subsection}


\begin{subsection}{More preparatory lemmas}

We need to show that the $\Z[\delta]$-module $\mathcal{M}[\Diag^{b}_{n}(\Omega)]$ is closed under multiplication, making it a $\Z[\delta]$-algebra.  First, we shall prove a few additional lemmas that will aid in the process.

\begin{figure}[!ht]
\centering
\subcaptionbox{Lemma~\ref{lem:closed under mult 1}.\label{Fig150}}{\includegraphics[scale=.95]{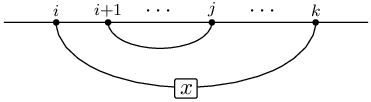}}
\qquad
\subcaptionbox{Lemma~\ref{lem:closed under mult 2}.\label{Fig151}}{\raisebox{9ex-\height}{\includegraphics[scale=.95]{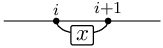}}}
\qquad
\subcaptionbox{Lemma~\ref{lem:closed under mult 3}.\label{Fig152}}{\raisebox{9.2ex-\height}{\includegraphics[scale=.95]{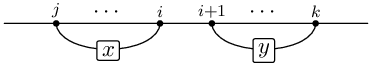}}}\\
\vspace{1ex}
\subcaptionbox{Lemma~\ref{lem:closed under mult 4}.\label{Fig154}}{\includegraphics[scale=.95]{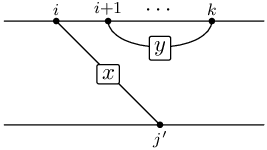}}
\qquad
\subcaptionbox{Lemma~\ref{lem:closed under mult 5}.\label{Fig155}}{\includegraphics[scale=.95]{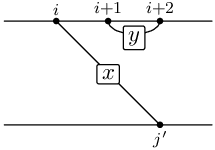}}
\qquad
\subcaptionbox{Lemma~\ref{lem:closed under mult 6}.\label{Fig161}}{\includegraphics[scale=.95]{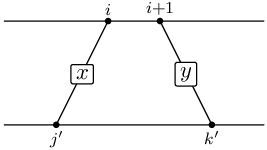}}
\caption{Edge configurations for Lemmas~\ref{lem:closed under mult 1}--\ref{lem:closed under mult 7}}\label{Fig150--Fig161}
\end{figure}

\begin{lemma}\label{lem:closed under mult 1}
Let $d$ be an admissible diagram with the edge configuration at nodes $i$ and $i+1$ given in Figure~\ref{Fig150}, where $x$ represents a (possibly trivial) block of decorations.  Then $d_{i}d=2^{c}d'$, where $c \in \{0,1\}$ and $d'$ is an admissible diagram.  Moreover, $c=1$ if and only if $i=1$.
\end{lemma}

\begin{proof}
The only case that requires serious consideration is if $i=1$; the result follows immediately if $i>1$.  Assume that $i=1$.  Since $d$ is admissible, $x \in \{\bcirc, \bcirc \wtri, \bcirc \wcirc\}$.  In any case, $d_{1}d=2 d'$ for some diagram $d'$, where the non-propagating edge joining node $j$ to node $k$ in $d'$ is one of the following blocks: $\btri, \btri \wtri$, or $\btri \wcirc$.  It follows that $d'$ is admissible.
\end{proof}

\begin{lemma}\label{lem:closed under mult 2}
Let $d$ be an admissible diagram with the edge configuration at nodes $i$ and $i+1$ given in Figure~\ref{Fig151}, where $x$ represents a (possibly trivial) block of decorations.  Then $d_{i}d=\delta^{c}d'$, where $c \in \{0,1\}$ and $d'$ is an admissible diagram.  Moreover, $c=0$ if and only if $x \in \{\bcirc \wtri, \btri \wtri, \btri \wcirc\}$.
\end{lemma}

\begin{proof}
We consider two cases.  For the first case, assume that $1<i<n+1$.  Since $d$ is admissible, $x \in \{\emptyset, \btri, \wtri, \btri \wtri\}$.  (Note that $x =\btri \wtri$ only if $d$ is undammed; otherwise $d$ would not be LR-decorated.)  In either case, $d_{i}d$ produces a loop decorated with the block $x$ along with a diagram that is identical to $d$, except that the block $x$ has been removed from the edge joining $i$ to $i+1$.  The loop decorated with the block $x$ is equal to $\delta$, unless $x=\btri \wtri$, in which case the loop is irreducible.  Regardless, the resulting diagram is admissible, as desired.  For the second case, assume that $i=1$ or $n$.  Without loss of generality, assume that $i=1$, the other case being symmetric.  Since $d$ is admissible, $x \in \{\bcirc, \bcirc \wtri\}$.  If $x=\bcirc$, then $d_{1}d=\delta d$, as expected.  If, on the other hand, $x=\bcirc \wtri$ (which can only happen if $d$ is undammed), then $d_{1}d$ results in an admissible diagram that is identical to $d$ except that we add a loop decorated by $\btri \wtri$ and remove the $\wtri$ decoration from the edge connecting node 1 to node 2.  
\end{proof}

\begin{lemma}\label{lem:closed under mult 3}
Let $d$ be an admissible diagram with the edge configuration at nodes $i$ and $i+1$ given in Figure~\ref{Fig152}, where $x$ and $y$ represent (possibly trivial) blocks of decorations.  Then $d_{i}d=2^{c}d'$, where $c \in \{0,1\}$ and $d'$ is an admissible diagram.  
\end{lemma}

\begin{proof}
First, observe that $d_{i}d$ has the edge configuration at nodes $i$ and $i+1$ given in Figure~\ref{Fig153}, where $xy=2^{c}z$ and $z$ is a basis element of $\V$.  Note that since $d$ is admissible, there will be at most one relation to apply in the product $xy$, which will happen exactly when the last decoration in $x$ and the first decoration in $y$ are of the same type (open or closed).  This implies that $c \in \{0,1\}$.  If $j=1$ (respectively, $k=n+2$), then the first (respectively, last) decoration in $x$ (respectively, $y$) must be a $\bcirc$ (respectively, $\wcirc$) decoration.  Furthermore, if $j=1$ (respectively, $k=n+2$), then this is the only occurrence of a $\bcirc$ (respectively, $\wcirc$) decoration on a non-propagating edge in the north face of $d$.  By inspecting the possible relations we can apply, this implies that if $j=1$ (respectively, $k=n+2$), the first (respectively, last) decoration of $z$ must be a $\bcirc$ (respectively, $\wcirc$) decoration and this is the only occurrence of a $\bcirc$ (respectively, $\wcirc$) decoration on a non-propagating edge of the diagram that results from the product $d_{i}d$.  If, on the other hand, $j\neq 1$ and $k\neq n+2$, then neither of $x$ or $y$ may contain a $\bcirc$ or $\wcirc$ decoration.  In this case, $z$ will not contain any $\bcirc$ or $\wcirc$ decorations either.  This argument shows that the diagram that results from the product $d_{i}d$ must be admissible.
\end{proof}

\begin{figure}[!ht]
\centering
$2^{c}\ \begin{tabular}[c]{@{}c@{}}
\includegraphics[scale=.95]{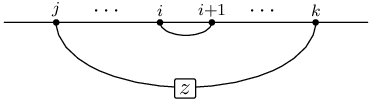}
\end{tabular}$
\caption{Diagram for the proof of Lemma~\ref{lem:closed under mult 3}.}\label{Fig153}
\end{figure}

\begin{lemma}\label{lem:closed under mult 4}
Let $d$ be an admissible diagram such that $\a(d)>1$ with the edge configuration at nodes $i$ and $i+1$ given in Figure~\ref{Fig154}, where $x$ and $y$ represent (possibly trivial) blocks of decorations.  Then $d_{i}d=2^{c}d'$, where $c \in \{0,1\}$ and $d'$ is an admissible diagram.  
\end{lemma}

\begin{proof}
Note that $1\leq i <n+1$.  Since $d$ is dammed, $y$ is either equal to the identity in $\V$ or is equal to an open decoration.  On the other hand, $x$ could be equal to the identity in $\V$, a single closed decoration, a single open decoration, or if $d$ has a unique propagating edge, then $x$ could be an alternating sequence of open and closed decorations.  We consider two cases: (1) $1<i<n+1$ and (2) $i=1$.

Case (1).  If $1<i<n+1$, then there will not be any relations to apply in the product of $d_{i}$ and $d$ unless the first decoration on the edge joining $i$ to $j'$ in $d$ is open and $y$ is also an open decoration.  In this case, $d_{i}d$ will be equal to 2 times an admissible diagram.  

Case (2).  Now, assume that $i=1$.  Since $d$ is admissible, either $x$ is trivial or the first decoration on the edge joining $1$ to $j'$ in $d$ must be closed.  If $x$ is trivial, then $j'=1'$, in which case, $d_{i}d$ is equal to a single admissible diagram.  If the first decoration is closed, then $d_{i}d$ equals 2 times an admissible diagram, as expected.
\end{proof}

\begin{lemma}\label{lem:closed under mult 5}
Let $d$ be an admissible diagram such that $\a(d)=1$ with the edge configuration at nodes $i$ and $i+1$ given in Figure~\ref{Fig155}, where $x$ and $y$ represent (possibly trivial) blocks of decorations.  Then $d_{i}d=2^{c}d'$, where $c \in \{0,1\}$ and $d'$ is an admissible diagram with $\a(d')=1$.  
\end{lemma}

\begin{proof}
Since $\a(d)=1$, the non-propagating edge joining $i+1$ to $i+2$ is the unique non-propagating edge in the north face of $d$. Furthermore, since $\a(d)=1$, the edge configuration at nodes $i$ and $i+1$ forces $j \in \{i, i+2\}$.  According to Lemma~\ref{lem:a-value=1}, the diagram that is produced by multiplying $d_{i}$ times $d$ has $\a$-value 1.  We consider three cases: (1) $i=1$, (2) $1< i<n$, and (3) $i=n$.

Case (1).  Assume that $i=1$.  This implies that $j \in \{1, 3\}$.  Then the possible edge configurations at nodes 1 and 2 of $d$ that are consistent with axiom~\ref{C3} of Definition~\ref{def:admissible} are the ones listed in Figures~\ref{Fig156} and~\ref{Fig157}, where the rectangle represents a (possibly trivial) sequence of blocks such that each block is a single $\btri$. In any case, we see that $d_{i}d=d_{1}d=2^{c}d'$, where $c \in \{0,1\}$ and $d'$ is an admissible diagram.

\begin{figure}[!ht]
\centering
\subcaptionbox{\label{Fig156}}{\includegraphics[scale=.95]{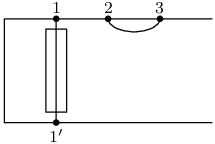}}
\qquad
\subcaptionbox{\label{Fig157}}{\includegraphics[scale=.95]{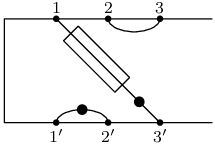}}\\
\subcaptionbox{\label{Fig158}}{\includegraphics[scale=.95]{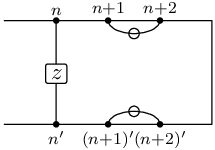}}
\qquad
\subcaptionbox{\label{Fig159}}{\includegraphics[scale=.95]{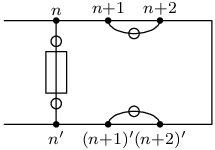}}
\qquad
\subcaptionbox{\label{Fig160}}{\includegraphics[scale=.95]{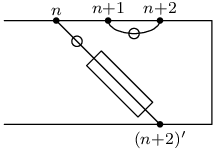}}
\caption{Diagrams for cases (1) and (3) of the proof of Lemma~\ref{lem:closed under mult 5}.}\label{Fig156--Fig160}
\end{figure}

Case (2).  Next, assume that $1<i<n$.  Since $\a(d)=1$, the restrictions on $i$ and $j'$ imply that both $x$ and $y$ are trivial.  That is, the propagating edge from $i$ to $j'$ and the non-propagating edge from $i+2$ to $i+3$ are undecorated.  Therefore, it is quickly seen that $d_{i}d=d'$ for some admissible diagram $d'$.

Case (3).  For the final case, assume that $i=n$.  This implies that $j \in \{n, n+2\}$, in which case the possible edge configurations at nodes $n$ and $n+1$ of $d$ that are consistent with axiom~\ref{C3} of Definition~\ref{def:admissible} are the ones listed in Figures~\ref{Fig158},~\ref{Fig159}, and~\ref{Fig160}, where $z\in \{\emptyset,\wtri\}$ and the rectangles on Figures~\ref{Fig159} and~\ref{Fig160} represent a sequence of blocks such that each block is a single $\wtri$.  In any case, we see that $d_{i}d=d_{n}d=2^{c}d'$, where $c \in \{0,1\}$ and $d'$ is an admissible diagram.
\end{proof}

\begin{lemma}\label{lem:closed under mult 6}
Let $d$ be an admissible diagram such that $\a(d)>1$ with the edge configuration at nodes $i$ and $i+1$ given in Figure~\ref{Fig161}, where $x$ and $y$ represent (possibly trivial) blocks of decorations.  Then $d_{i}d=2^{c}d'$, where $c \in \{0,1\}$ and $d'$ is an admissible diagram.  
\end{lemma}

\begin{proof}
Since $d$ is LR-decorated, $x$ and $y$ cannot be of the same type (open or closed).  The only time there is potential to apply any relations when multiplying $d_{i}$ times $d$ is if $i=1$ (respectively, $i=n+1$) and $x$ (respectively, $y$) is nontrivial.  Regardless, it is easily seen that the statement of the lemma is true.
\end{proof}

\begin{lemma}\label{lem:closed under mult 7}
Let $d$ be an admissible diagram such that $\a(d)=1$ with the edge configuration at nodes $i$ and $i+1$ given in Figure~\ref{Fig161}, where $x$ and $y$ represent sequences of (possibly trivial) blocks of decorations.  Then $d_{i}d=2^{k}d'$, where $k\geq 0$ and $d'$ is an admissible diagram with $\a(d)>1$. 
\end{lemma}

\begin{proof}
According to Lemma~\ref{lem:a-value=1}, the diagram that is produced by multiplying $d_{i}$ times $d$ has $\a$-value strictly greater than 1.  In this case, the sequence of blocks of decorations occurring on the leftmost (respectively, rightmost) propagating edge of $d$ will conjoin in the product of $d_{i}$ and $d$.  This implies that $d_{i}d=2^{k}d'$ for $k \geq 0$ and some diagram $d'$.  To see that $d'$ is admissible, we consider the five possibilities for $d$ given in Figure~\ref{Fig163--Fig167}, where $u\in \{\emptyset,\btri\}$ and the rectangle on the leftmost (respectively, rightmost) propagating edge represents a (possibly trivial) sequence of blocks such that each block is a single $\btri$ (respectively, $\wtri$); all remaining possibilities are analogous.

In each of these cases, if $d$ has propagating edges joined to nodes $i$ and $i+1$ in the north face,  it is quickly seen that the diagram $d'$ that results from multiplying $d_{i}$ times $d$ will be consistent with the axioms of Definition~\ref{def:admissible} since $\bcirc \btri \cdots \btri \bcirc$ and $\btri \cdots \btri$ (respectively, $\wcirc \wtri \cdots \wtri \wcirc$ and $\wtri \cdots \wtri$) are equal to a power of 2 times $\btri$ (respectively, $\wtri$).
\end{proof}

\end{subsection}


\begin{subsection}{The admissible diagrams form a basis}

The next proposition states that the product of a simple diagram and an admissible diagram results in a multiple of an admissible diagram.  The proof relies on stringing together Lemmas~~\ref{lem:closed under mult 1}--\ref{lem:closed under mult 7}.

\begin{proposition}\label{prop:powers of 2 and delta}
Let $d$ be an admissible diagram.  Then $d_{i}d=2^{k}\delta^{m}d'$ for some $k,m \in \Z^{+}\cup\{0\}$ and admissible diagram $d'$.
\end{proposition}

\begin{proof}
Let $d$ be an admissible diagram and consider the product $d_{i}d$.  Observe that the only possible edge configurations for $d$ at nodes $i$ and $i+1$ are the ones in Figure~\ref{Fig168--Fig174}.

\begin{figure}[!ht]
\centering
\subcaptionbox{\label{Fig168}}{\includegraphics[scale=.95]{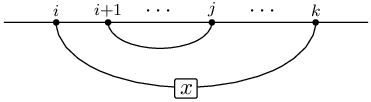}}
\qquad
\subcaptionbox{\label{Fig169}}{\includegraphics[scale=.95]{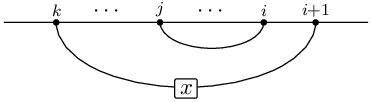}}\\
\subcaptionbox{\label{Fig170}}{\raisebox{5.5ex-\height}{\includegraphics[scale=.95]{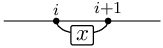}}}
\qquad
\subcaptionbox{\label{Fig171}}{\includegraphics[scale=.95]{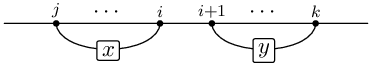}}\\
\subcaptionbox{\label{Fig172}}{\includegraphics[scale=.95]{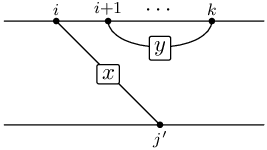}}
\qquad
\subcaptionbox{\label{Fig173}}{\includegraphics[scale=.95]{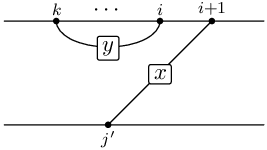}}
\qquad
\subcaptionbox{\label{Fig174}}{\includegraphics[scale=.95]{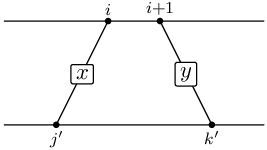}}
\caption{The seven possible edge configurations of Proposition~\ref{prop:powers of 2 and delta}.}\label{Fig168--Fig174}
\end{figure}

If $d$ is the diagram in Figure~\ref{Fig168}, the result follows from Lemma~\ref{lem:closed under mult 1}, and if $d$ is the diagram in Figure~\ref{Fig169}, we may apply a symmetric argument.  In the case of Figure~\ref{Fig170}, the result follows from Lemma~\ref{lem:closed under mult 2}.  Lemma~\ref{lem:closed under mult 3} may be applied when $d$ is the diagram in Figure~\ref{Fig171}.  If $d$ is the diagram in Figure~\ref{Fig172}, we need only apply Lemmas~\ref{lem:closed under mult 4} and~\ref{lem:closed under mult 5}, and when $d$ is the diagram in Figure~\ref{Fig173} the result follows by a symmetric argument.  Finally, Lemmas~~\ref{lem:closed under mult 6} and~\ref{lem:closed under mult 7} handle the case when $d$ is the diagram in Figure~\ref{Fig174}.
\end{proof}

\begin{corollary}\label{cor:module is subalgebra}
The $\Z[\delta]$-module $\mathcal{M}[\Diag^{b}_{n}(\Omega)]$ is a $\Z[\delta]$-subalgebra of $\widehat{\P}_{n+2}^{LR}(\Omega)$.
\end{corollary}

\begin{proof}
This statement follows immediately from Propositions~~\ref{prop:admissibles gen by simples} and~\ref{prop:powers of 2 and delta}.
\end{proof}

We are finally ready to show that the admissible diagrams form a basis for $\D_n$.

\begin{theorem}\label{thm:module admissibles is D_n}
The $\Z[\delta]$-algebras $\mathcal{M}[\Diag^{b}_{n}(\Omega)]$ and $\D_{n}$ are equal.  Moreover, the set of admissible diagrams is a basis for $\D_{n}$.
\end{theorem}

\begin{proof}
Proposition~\ref{prop:admissibles gen by simples} and Corollary~\ref{cor:module is subalgebra} imply that $\mathcal{M}[\Diag^{b}_{n}(\Omega)]$ is a subalgebra of $\D_{n}$.  However, $\D_{n}$ is the smallest algebra containing the simple diagrams, which $\mathcal{M}[\Diag^{b}_{n}(\Omega)]$ also contains since the simple diagrams are admissible.  Therefore, we must have equality of the two algebras.  By Proposition~\ref{prop:admissible basis}, the set of admissible diagrams is a basis for $\mathcal{M}[\Diag^{b}_{n}(\Omega)]$.  Therefore, the set of admissible diagrams forms a basis for $\D_{n}$.
\end{proof}

\end{subsection}

\end{section}


\begin{section}{Closing remarks}\label{sec:closing}

In this paper, we constructed an infinite dimensional associative diagram algebra $\D_n$.  We were able to easily check that this algebra satisfies the relations of $\TL(\C_n)$, thus showing that there is a surjective algebra homomorphism from $\TL(\C_n)$ to $\D_n$.  Moreover, we described the set of admissible diagrams and accomplished the more difficult task of proving that this set of diagrams forms a basis for $\D_n$.  

What remains to be shown is that our diagrammatic representation is faithful and that each admissible diagram corresponds to a unique monomial basis element of $\TL(\C_n)$.  Demonstrating injectivity of the homomorphism between $\TL(\C_n)$ and $\D_n$ is dealt with in the sequel to this paper~\cite{Ernst.D:D} (also see~\cite{Ernst.D:A}).

One motivation behind studying these generalized Temperley--Lieb algebras is that they provide a gateway to understanding the Kazhdan--Lusztig theory of the associated Hecke algebra.  Recall that if $(W,S)$ is Coxeter system of type $\Gamma$, the associated Hecke algebra $\H(\Gamma)$ is an algebra with a basis given by $\{T_w: w \in W\}$ and relations that deform the relations of $W$ by a parameter $q$. Loosely speaking, $\TL(\Gamma)$ retains some of the relevant structure of $\H(\Gamma)$, yet is small enough that computation of the leading coefficients of the notoriously difficult to compute Kazhdan--Lusztig polynomials is often much simpler.

Using the diagrammatic representations of $\TL(\Gamma)$ when $\Gamma$ is of types $A, B, D,$ or $E$, Green has constructed a trace on $\H(\Gamma)$ similar to Jones' trace in the type $A$ situation~\cite{Green.R:K,Green.R:L}.  Remarkably, this trace can be used to non-recursively compute leading coefficients of Kazhdan--Lusztig polynomials indexed by pairs of FC elements, and this is precisely our motivation in the type $\C$ case.  

In a future paper, we plan to construct a Jones-type trace on $\H(\C)$ using the diagrammatic representation of $\TL(\C)$, thus allowing us to be able to quickly compute leading coefficients of the infinitely many Kazhdan--Lusztig polynomials indexed by pairs of FC elements.  Understanding the diagrammatic representation of $\TL(\C_n)$ and its corresponding Jones-type trace should provide insight into what happens in the more general case involving an arbitrary Coxeter graph $\Gamma$.

\end{section}


\section*{Acknowledgements}

I would like to thank R.M.~Green for many useful conversations during the preparation of this article.  I am also grateful to the referee for his or her careful reading of the paper and constructive suggestions for improvements. 

\bibliographystyle{plain}
\bibliography{DiagCalc1}

\end{document}